\documentclass[12pt]{amsart}
\usepackage{amssymb}
\usepackage{amsthm}
\usepackage{bm}
\usepackage{latexsym}
\usepackage{amsmath}
\usepackage{eufrak}
\usepackage{mathrsfs}
\usepackage{amscd}
\usepackage[all]{xy}
\usepackage[usenames]{color}
\usepackage[dvips]{graphicx}

\theoremstyle{plain}
\newtheorem{teo}{Theorem}[section]
\newtheorem{cor}[teo]{Corollary}
\newtheorem{lem}[teo]{Lemma}
\newtheorem{question}[teo]{Question}
\newtheorem{prop}[teo]{Proposition}
\theoremstyle{definition}
\newtheorem{defi}{Definition}[section]
\newtheorem{obs}{Remark}[section]
\newtheorem{example}{Example}[section]

\DeclareMathOperator{\aut}{{Aut}}

\DeclareMathOperator{\spec}{Spec}
\DeclareMathOperator{\sing}{sing}

\def\B{{\mathbb{B}}}
\def\FF{{\mathcal{F}}}

\def\O{{\mathcal{O}}}

\def\L{{\mathcal{L}}}
\def\F{{\mathbb{F}}}
\def\K{{\mathbb{K}}}

\def\M{{\mathcal{M}}}
\def\N{{\mathcal{N}}}
\def\A{{\mathcal{A}}}

\def\AA{{\mathscr{A}}}

\def\LL{{\mathscr{L}}}

\def\E{{\mathcal E}}
\def\Z{{\mathbb Z}}
\def\Q{{\mathbb Q}}
\def\H{{\mathscr H}}

\def\P{{\mathbb P}}
\def\C{{\mathbb C}}
\def\Af{{\mathbb A}}
\def\EE{{\mathbb{E}}}

\def\spec{\textrm{Spec}}
\def\Pic{\textrm{Pic}}

\begin{document}
\bibliographystyle{amsplain}

\title{Arrangements of rational sections over curves and the varieties they define}
\author{\textrm{Giancarlo Urz\'ua}}
\date{\today}

\email{urzua@math.umass.edu} \maketitle

\begin{abstract}
We introduce arrangements of rational sections over curves. They
generalize line arrangements on $\P^2$. Each arrangement of $d$
sections defines a single curve in $\P^{d-2}$ through the
Kapranov's construction of $\overline{M}_{0,d+1}$. We show a
one-to-one correspondence between arrangements of $d$ sections and
irreducible curves in $M_{0,d+1}$, giving also correspondences for
two distinguished subclasses: transversal and simple crossing.
Then, we associate to each arrangement $\A$ (and so to each
irreducible curve in $M_{0,d+1}$) several families of nonsingular
projective surfaces $X$ of general type with Chern numbers
asymptotically proportional to various log Chern numbers defined
by $\A$. For example, for the main families and over $\C$, any
such $X$ is of positive index and $\pi_1(X) \simeq
\pi_1(\overline{A})$, where $\overline{A}$ is the normalization of
$A$. In this way, any rational curve in $M_{0,d+1}$ produces
simply connected surfaces with $2< \frac{c_1^2(X)}{c_2(X)} <3$.
Inequalities like these come from log Chern inequalities, which
are in general connected to geometric height inequalities (see
Appendix). Along the way, we show examples of \'etale simply
connected surfaces of general type in any characteristic violating
any sort of Miyaoka-Yau inequality.

\end{abstract}

\vspace{0.5cm}

Arrangements of rational sections over curves set up a new class
of arrangements of curves on algebraic surfaces. Given a
nonsingular projective curve $C$ and an invertible sheaf $\L$ on
$C$, they are defined as finite collections of sections of the
corresponding $\Af^1$-bundle. The simplest example is line
arrangements on $\P^2$, where $C=\P^1$ and $\L=\O_{\P^1}(1)$. In
Section \ref{s1}, we systematically study arrangements of rational
sections over curves. Although in somehow they can be managed as
line arrangements, the big difference relies on possible
tangencies among their curves, introducing more geometric
liberties. We partially organize this via transversal and simple
crossing arrangements (Definition \ref{d14}). Throughout Sections
\ref{s2}, \ref{s3} and \ref{s4}, we show one-to-one
correspondences between arrangements of $d$ sections and
irreducible curves in $M_{0,d+1}$, the moduli space of curves of
genus zero with $d+1$ marked ordered points. This is done for each
fixed pair $(C,\L)$ in the general (Theorem \ref{t41}),
transversal (Corollary \ref{c44}), and simple crossing (Corollary
\ref{c45}) cases. Because of Kapranov's description of
$\overline{M}_{0,d+1}$ \cite{Ka1,Ka2}, this produces bijections
between arrangements and curves in $\P^{d-2}$ (Corollary
\ref{c41}). For instance, arrangements of $d$ lines in $\P^2$
correspond to lines in $\P^{d-2}$ (as in \cite{U1}), arrangements
of $d$ conics in $x^2+y^2+z^2=0$ correspond to conics in
$\P^{d-2}$, etc. To exemplify, we show in Section \ref{s5} a way
to produce explicit arrangements of sections from irreducible
curves in $\P^2$. This is based on \cite[Section 7]{CT1}, where
the authors show how to cover $M_{0,d+1}$ with blow-ups of $\P^2$
at $d+1$ points. We use their rigid conic as concrete example (see
Examples \ref{e51} and \ref{e61}).

Given an arrangement of sections $\A$, we define two types of
arrangements: the extended $\A_{\Delta}$ and some partially
extended $\A_{p\Delta}$. Their definitions and log properties are
exposed in Section \ref{s6}. Over $\C$, they satisfy certain log
Miyaoka-Yau inequalities which are no longer combinatorial as in
the case of line arrangements (Remarks \ref{o64} and \ref{o62}).
For line arrangements we know an optimal log inequality
(Hirzebruch-Sakai's in Remark \ref{o64}), but for any other class
we only have the coarse bound $3$. Arrangements attaining upper
bounds should be special, and would produce interesting algebraic
surfaces by means of Theorem \ref{t71}. We remark that questions
about sharp upper bounds of log Chern ratios are related to
questions on effective height inequalities
\cite[pp.149-153]{Lang91} (see Appendix, where we slightly extend
and give another proof of Liu's inequality \cite[Theorem
0.1]{Liu96}, which naturally implies strict Tan's height
inequality \cite[Theorem A]{Tan95}).

In Section \ref{s7}, we associate various families of nonsingular
projective surfaces to any given arrangement of sections $\A$.
These surfaces share the random nature of the surfaces $X$
constructed in \cite{U2}, having Chern invariants asymptotically
proportional to the log Chern invariants of $\A_{\Delta}$ and
$\A_{p\Delta}$'s. In this way, we are able to show for a more
general class of arrangements (and so singularities) that the
behavior of Dedekind sums and continued fractions used in
\cite{U2} can also be applied. In this paper, any such $X$ is of
general type and satisfies $c_1^2(X),c_2(X)>0$. Putting it all
together, and over $\C$, we have the following.

\vspace{0.3 cm} \textbf{Theorem.} Let $A$ be an irreducible curve
in $\P^n$ not contained in $\prod_i x_i \prod_{i<j} (x_j-x_i)=0$.
Let $\overline{A}$ be the normalization of $A$. Then, there exist
nonsingular projective surfaces $X$ of general type such that $2 <
\frac{c_1^2(X)}{c_2(X)} < 3,$ having $\frac{c_1^2(X)}{c_2(X)}$
arbitrarily close to $\frac{\bar{c}_1^2(A)}{\bar{c}_2(A)}$, a
well-defined positive rational number depending on $A$ and its
position in $\P^n$. In addition, there is an induced connected
fibration $\pi': X \rightarrow \overline{A}$ which gives an
isomorphism $\pi_1(X) \simeq \pi_1(\overline{A})$. In this way,
Alb$(X) \simeq$ Jac$(\overline{A})$ and $\pi'$ is the Albanese
fibration of $X$. \vspace{0.3 cm}

With this in hand, we aim to answer the still open question: are
there simply connected nonsingular projective surfaces of general
type with $\frac{c_1^2}{c_2}$ arbitrarily close to the Miyaoka-Yau
bound $3$? Hence, at least when $A$ is a rational curve, it is
important for us to know about sharp upper bounds for
$\frac{\bar{c}_1^2(A)}{\bar{c}_2(A)}$ (also for $\A_{p\Delta}$,
see Remark \ref{o73}). So far, we only know that this bound is
strictly smaller than $3$ (Corollaries \ref{c61} and \ref{c62}).
On the other hand, in positive characteristic, we use our method
to produce \'etale simply connected nonsingular projective
surfaces of general type which violate any sort of Miyaoka-Yau
inequality for any given characteristic (Example \ref{e15}, Remark
\ref{o65}, Example \ref{e72}).

\vspace{0.3cm}

\textbf{Acknowledgments}: I am grateful to I. Dolgachev, S. Keel,
F. Knudsen, K. Liu, N. Mok, A. Momot, J. Tevelev, T. Weston and S.
Wong for valuable conversations and e-mail correspondence. I would
like to thank the referee for such a complete and helpful review
which improved much many parts of this article.   

\tableofcontents

\section{Arrangements of rational sections over curves} \label{s1}

Fix an algebraically closed field $\K$. Let $C$ be a nonsingular
projective curve defined over $\K$ of genus $g=h^1(C,\O_C)$.
Hence, when $\K=\C$, the curve $C$ is a compact Riemann surface of
genus $g$. Let $\L$ be an invertible sheaf on $C$ of degree
$\deg(\L)=e>0$, and let
$$\Af_C(\L):= \spec(S(\L^{-1})) \rightarrow C$$ be the line bundle associated to $\L$ (as in \cite[II,Ex.5.18]{Ha1}), where $\L^{-1}$ is the dual
of $\L$. This is the so-called total space of $\L$. A section of
$\Af_C(\L) \rightarrow C$ is a morphism $C \rightarrow \Af_C(\L)$
such that the composition map $C \rightarrow \Af_C(\L) \rightarrow
C$ is the identity. The space of sections can be identified with
$H^0(C,\L)$. Since it is better to deal with a projective surface,
we naturally compactify all fibers, so that we work with a
$\P^1$-bundle. Let $$\pi: \P(\O_C \oplus \L^{-1}) \rightarrow C$$
be the $\P^1$-bundle associated to $\O_C \oplus \L^{-1}$ over $C$,
as in \cite[II,Ex.7.10]{Ha1}. The nonsingular projective surface
$\P(\O_C \oplus \L^{-1})$ contains $\Af_C(\L)$ as an open set,
such that the curve $C_0 := \P(\O_C \oplus \L^{-1}) \setminus
\Af_C(\L)$ is a section of $\pi$ with self-intersection
$C_0^2=-e$. It is easy to see that $C_0$ is the only irreducible
curve with negative self-intersection in $\P(\O_C \oplus
\L^{-1})$. This surface is a particular case of a geometrically
ruled surface over $C$ \cite[V,Section 2]{Ha1}, and it is in its
normalized form \cite[V,Proposition 2.8]{Ha1}. We denote by $F_c$
the fiber over a point $c \in C$, or just $F$ when we consider its
numerical class. Any element in $\Pic \big( \P(\O_C \oplus
\L^{-1}) \big)$ can be written as $aC_0 + \pi^*(\M)$ with $a \in
\Z$, and $\M \in \Pic(C)$. Any element in $\text{Num} \big(
\P(\O_C \oplus \L^{-1}) \big)$ can be written as $aC_0 + bF$ with
$a,b \in \Z$ \cite[p.373]{Ha1}.

\begin{example}
Let $C=\P^1$, and let $e>0$. Consider the invertible sheaf
$\O_{\P^1}(e)$ on $\P^1$. Then, the surface $\P \big(\O_{\P^1}
\oplus \O_{\P^1}(-e)\big)$ is the Hirzebruch surface $\F_e$. When
$e=1$, $\F_1$ corresponds to the blow-up at a point of $\P^2$
\cite[V,Exa.2.11.4]{Ha1}, and $C_0$ is the $(-1)$-curve.
\label{e11}
\end{example}

The main objects are the following.

\begin{defi}
Let $d\geq 3$ be an integer. Let $C$ be a nonsingular projective
curve, and let $\L$ be an invertible sheaf on $C$ of degree $e>0$.
An \textbf{arrangement of $d$ sections} is a labeled set of $d$
distinct sections $\A= \{ S_1, \ldots, S_d \}$ of $\pi: \P(\O_C
\oplus \L^{-1}) \rightarrow C$ such that $$S_i \sim C_0 +
\pi^*(\L)$$ for all $i \in \{1,2,\ldots ,d \}$, and
$\bigcap_{i=1}^d S_i = \emptyset$. From now on, we denote $\P(\O_C
\oplus \L^{-1})$ by $\P_C(\L)$. \label{d11}
\end{defi}

In particular, $S_i.S_j=e$, and $S_i.C_0=0$ for all $i$, and so
these arrangements are indeed formed by sections of $\Af_C(\L)
\rightarrow C$. The condition $\bigcap_{i=1}^d S_i = \emptyset$
implies that $\L$ is base-point free. To see this, take a point $c
\in C$, and consider the corresponding fiber $F_c$. Since
$\bigcap_{i=1}^d S_i = \emptyset$, there are two sections $S_i$,
$S_j$ such that $F_c \cap S_i \cap S_j = \emptyset$. Let
$\sigma_j: C \rightarrow \P_C(\L)$ be the map defining the section
$S_j$. Then, $\L \simeq \sigma_j^*(\pi^*(\L)\otimes \O_{S_j})
\simeq \sigma_j^*(\O_{\P_C(\L)}(C_0) \otimes \pi^*(\L)\otimes
\O_{S_j}) \simeq \sigma_j^*(\O_{\P_C(\L)}(S_i) \otimes \O_{S_j})$,
and $\sigma_j^*(\O_{\P_C(\L)}(S_i) \otimes \O_{S_j})$ is given by
an effective divisor on $C$ not supported at $c$. This tell us
that $\L \simeq \O_C(D)$ with $D$ base-point free effective
divisor.

If $\A = \{S_1, \ldots, S_d \}$ is an arrangement as in Definition
\ref{d11}, but with $\bigcap_{i=1}^d S_i \neq \emptyset$, then we
can apply elementary transformations (see \cite[V,Exa.5.7.1]{Ha1})
at each of the points in $\bigcap_{i=1}^d S_i$ to obtain a new
arrangement of $d$ sections in $\P_C(\L')$ for some $\L'$. After
repeating this process a finite number of times, we arrive to an
arrangement $\A'$ in $\P_C(\L')$ with $\bigcap_{i=1}^d S'_i =
\emptyset$. If $\deg(\L')=0$, then $\L'=\O_C$ since, as we showed
above, $\L'=\O_C(D)$ for some effective divisor $D$. In this case,
$\P_C(\O_C) = C \times \P^1$, and the arrangement is trivially
formed by a collection of $d$ ``horizontal" fibers (it just
corresponds to an arrangement of $d$ points in $\P^1$). If $\A =
\{S_1, \ldots, S_d \}$ is a collection of arbitrary $d$ sections
in $\P_C(\L)$, we perform elementary transformations on the points
in $C_0 \cap S_i$ for all $i$, and we repeat this process until
all sections are disjoint from the new curve $C'_0$ in $\P_C(\L')$
(proper transform of $C_0$). In this way, \textbf{arbitrary
arrangements of sections} can always be considered, after some
elementary transformations, as the ones in Definition \ref{d11}.

We now define the morphisms between our objects.

\begin{defi}
Fix an integer $d\geq 3$. Let $C,C'$ be nonsingular projective
curves, and let $\L, \L'$ be invertible sheaves of positive
degrees on $C,C'$ respectively. Let $\A,\A'$ be arrangements of
$d$ sections in $\P_C(\L), \P_{C'}(\L')$ respectively. A
\textbf{morphism of arrangements} is the existence of a finite map
$g: C \rightarrow C'$, and a commutative diagram
\begin{center} $ \xymatrix{ \P_C(\L) \ar[d]_{\pi} \ar[r]^{G} & \P_{C'}(\L') \ar[d]_{\pi'}   \\ C
\ar[r]^{g} & C' }$
\end{center} so that $\P_C(\L)$ is isomorphic to the base change by $g$, and $S_i=G^*(S'_i)$ for all $i$. If $g$ is
an isomorphism, then the arrangements are said to be isomorphic.
\label{d12}
\end{defi}

In particular, a curve $C$ with an automorphism $g$ produces
isomorphic arrangements via the pull back of $g$.

\begin{lem}
A morphism of arrangements satisfies $C_0=G^*(C'_0)$ and $g^*(\L')
\simeq \L$. We have $C_0^2=\deg(g) {C'}_0^2$ and $S_i^2=\deg(g)
{S'}_i^2$ for all $i$. \label{l11}
\end{lem}

\begin{proof}
Since $0=G^*(C'_0).G^*(S'_i)=G^*(C'_0).S_i$, we have that
$G^*(C'_0)=C_0$. We know that $\pi_*(C_0)= \O_C \oplus \L^{-1}$
and $\pi'_*(C'_0)= \O_{C'} \oplus \L'^{-1}$ (see
\cite[II,Proposition 7.11]{Ha1}). By flat base change
\cite[III,Proposition 9.3]{Ha1}, we have $$ g^* \pi'_*(C'_0)
\simeq \pi_* G^*(C'_0),$$ and so $g^*(\L') \simeq \L$. Therefore,
$\deg(\L)=\deg(g) \deg(\L')$, and so $C_0^2=\deg(g){C'}_0^2$, and
$S_i^2=\deg(g){S'}_i^2$ for all $i$.
\end{proof}

One wants to consider arrangements of sections which do not come
from others via base change, and so the following definition.

\begin{defi}
Let us fix the data $(C,\L,d)$ as above. An arrangement of $d$
sections $\A$ is said to be \textbf{primitive} if whenever we have
an arrangement $\A'$ for some data $(C',\L',d)$, and a morphism
$g$ as in Definition \ref{d12}, then $g$ is an isomorphism. The
set of \textbf{isomorphism classes} of primitive arrangements is
denoted by $\AA(C,\L,d)$. This is clearly independent of the
isomorphism class of $C$ and $\L$. \label{d13}
\end{defi}

For instance, if $\L$ has a base-point, then $\AA(C,\L,d)=
\emptyset$.

\begin{example}
Let $d\geq 3$ be an integer. An \textbf{arrangement of $d$ lines}
in the plane is a set of $d$ labeled lines $\A=\{L_1, \ldots, L_d
\}$ in $\P^2$ such that $\bigcap_{i=1}^d L_i = \emptyset$. As in
\cite{U1}, we introduced ordered pairs $(\A,P)$, where $\A$ is an
arrangement of $d$ lines, and $P$ is a point in $\P^2 \setminus
\bigcup_{i=1}^d L_i$. If $(\A,P)$ and $(\A',P')$ are two such
pairs, we say that they are isomorphic if there exists an
automorphism $T$ of $\P^2$ such that $T(L_i)=L'_i$ for every $i$,
and $T(P)=P'$. Given $(\A,P)$, we blow up the point $P$ to obtain
an arrangement of $d$ sections for the data $(\P^1,\O(1),d)$, and
given such an arrangement of sections, we blow down $C_0$ to get a
pair $(\A,P)$, where $P$ is the image of $C_0$ and $\A$ is formed
by the images of the sections. One sees that the set of pairs up
to isomorphism of pairs is precisely $\AA(\P^1,\O(1),d)$. By Lemma
\ref{l11}, any arrangement of $\P_{\P^1}(\O(1)) \rightarrow \P^1$
is primitive (degree considerations). This is the simplest case
for arrangements of rational sections over curves. Notice that
$$\AA(\P^1,\O(1),3) = \big(\{x_1=0,x_2=0,x_3=0\},[1:1:1] \big),$$ where $[x_1:x_2:x_3]$ are coordinates for $\P^2$.
\label{e12}
\end{example}

In the next sections, we will classify all primitive arrangements,
and some distinguished subclasses which are defined through
intersection properties of their members. We now look at these
intersections. In what follows, until the end of this section, we
fix the data $(C,\L,d)$.

\begin{defi}
Let $\A=\{ S_1, \ldots S_d \}$ be an arrangement of sections in
$\P_C(\L)$. Let $P$ be a point in $\P_C(\L)$, and let $f,g$ be
local equations defining $S_i,S_j$ at $P$. As in \cite[V,Section
1]{Ha1}, we define the \textbf{intersection multiplicity}
$(S_i.S_j)_P$ of $S_i$ and $S_j$ at $P$ to be the length of
$\O_{P,\P_C(\L)}/(f,g)$. If $P$ is not in $S_i$ or $S_j$, then
$(S_i.S_j)_P=0$. Notice that, since $S_i.S_j=e$, we have $0 \leq
(S_i.S_j)_P \leq e$. We distinguish two classes of arrangements:

\begin{itemize}
\item[(t)] We say that $\A$ is \textbf{transversal} if for any
$i\neq j$ and any point $P$ in $S_i \cap S_j$, there is $k\neq
i,j$ such that $(S_i.S_k)_P=(S_i.S_j)_P - 1$. The set of
isomorphism classes of primitive transversal arrangements is
denoted by $\AA_t(C,\L,d)$.

\item[(s)] We say that $\A$ is \textbf{simple crossing}
\footnote{these are the type of singularities for arrangements in
\cite{U2}.} if for any $i\neq j$ and any point $P$ in $S_i \cap
S_j$, we have $(S_i.S_j)=1$. This is, the members of the
arrangement are pairwise transversal. The set of isomorphism
classes of primitive simple crossing arrangements is denoted by
$\AA_s(C,\L,d)$.
\end{itemize}
\label{d14}
\end{defi}

\begin{obs}
In $(t)$ above, we have the requirement $(S_i.S_k)_P=(S_i.S_j)_P -
1$. This implies $(S_i.S_k)_P=(S_j.S_k)_P$, and so the definition
is symmetric on $i,j$. To see this, let $\sigma: Bl_P(\P_C(\L))
\rightarrow \P_C(\L)$ be the blow-up at $P$. Let $\widetilde{S}_a$
be the strict transforms of $S_a$, so that $\widetilde{S}_a \sim
\sigma^*(S_a) - E$, for $a=i,j,k$. Here $E$ is the exceptional
curve of $\sigma$. In this way, we have
$$\widetilde{S}_a.\widetilde{S}_b= S_a.S_b -1$$ since $S_a$ is
nonsingular at $P$. Since $\sigma$ is an isomorphism outside of
$E$, we have that
$(\widetilde{S}_a.\widetilde{S}_b)_{\widetilde{P}}=(S_a.S_b)_P
-1$, where $\widetilde{P}=\widetilde{S}_a \cap E$. If
$(S_i.S_j)_P=2$, then $(S_i.S_k)_P=1$, and so $(S_j.S_k)_P=1$. One
proves the general assertion by induction on $(S_i.S_j)_P$.
\label{r11}
\end{obs}

This gives the stratification $$\AA_s(C,\L,d) \subseteq
\AA_t(C,\L,d) \subseteq \AA(C,\L,d).$$ Notice that for line
arrangements $\AA_s(\P^1,\O(1),d) = \AA_t(\P^1,\O(1),d) =
\AA(\P^1,\O(1),d)$, but already for $(\P^1,\O(2),d)$ we have
different sets, as the next example shows.

\begin{example}
Consider collections of curves in $\P^2$ given by $A_i = \{ C_1 ,
C_2 , C_3 , C_4 \}$, as shown in Figure \ref{f1}. Here, $C_1$ is a
conic, and $C_2,C_3,C_4$ are lines. For distinct $i$'s, we have
different intersections among $C_j$'s. Each $A_i$ has a marked
point $P$ in $C_1$. Out of these configurations, we produce three
arrangements of sections in $\F_2$. We blow up $P$, and then we
perform an elementary transformation at $\widetilde{P}$, which is
the intersection of the strict transform of $C_1$ with the
exceptional divisor $E$. Then, we have an arrangement of sections
$\A_i = \{ S_1,S_2,S_3,S_4 \}$ in $\F_2$, where $S_j$ corresponds
to the strict transform of $C_j$.

\begin{figure}[htbp]
\includegraphics[width=11cm]{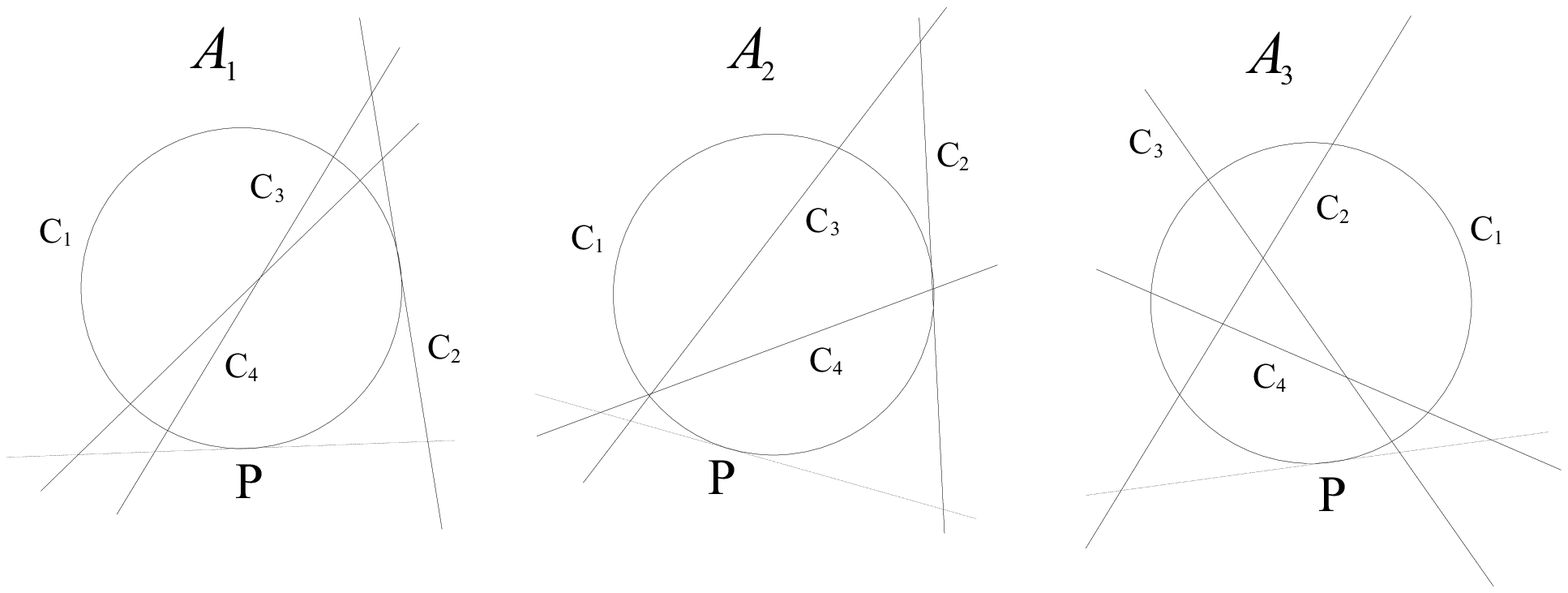} \caption{Configurations of curves in $\P^2$ which produce arrangements in $\F_2$.} \label{f1}
\end{figure}

Any possible morphism of arrangements, from $\A_i$ to some
$\A'_i$, would have $(\P^1,\O(1),4)$ as target, and the degree of
$g$ would have to be $2$. For $\A_1$, we have $8$ points in $\F_2$
where exactly two sections intersect, and $1$ where exactly $3$
intersect, so $\A_1$ is impossible as pull-back of $4$ sections in
$\P_{\P^1}(\O(1))$. Similar arguments apply to $\A_2$ and $\A_3$,
and so one easily checks that all of them a primitive. Notice that
$\A_3$ is simple crossing, $\A_2$ is a transversal, and $\A_1$ is
neither, so $\AA_s(\P^1,\O(2),4) \varsubsetneq \AA_t(\P^1,\O(2),4)
\varsubsetneq \AA(\P^1,\O(2),4)$. \label{e14}
\end{example}

\begin{figure}[htbp]
\includegraphics[width=10cm]{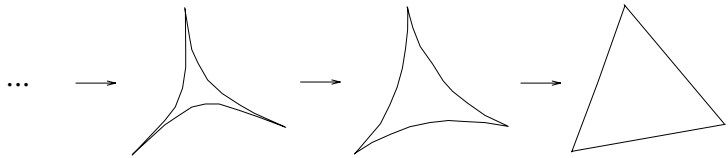} \caption{Evolution of a triangle under Frobenius maps in Example \ref{e15}.}
\end{figure}

\begin{example}
Assume $\K$ has characteristic $p>0$. Let $\A' \in
\AA_s(C',\L',d)$, and consider the $\K$-linear Frobenius morphism
$g:C \rightarrow C'$ \cite[p.301]{Ha1}, so $C$ and $C'$ are
isomorphic as abstract curves. Let $\A \in \AA(C,\L,d)$ be the
pull-back arrangement by $g$, as in Definition \ref{d12}. Then,
$g^*(\L')=\L$, and for any two members $S_i,S_j$ we have
$(S_i.S_j)_P=p$ when $P \in S_i \cap S_j$. The simple crossing
arrangement $\A'$ is transformed into an arrangement $\A$ where
any two members are tangent at $e$ points, each of order $p$.
\label{e15}
\end{example}

\section{Some facts about $\overline{M}_{0,d+1}$} \label{s2}

Let $d\geq 3$ be a integer. We denote by $\overline{M}_{0,d+1}$
the moduli space of $(d+1)$-pointed stable curves of genus zero
\cite{Knudsen,Ka1}. This is a smooth rational projective variety
of dimension $d-2$. The open set $M_{0,d+1}$ parametrizes
configurations of $d+1$ distinct labeled points in $\P^1$. The
boundary $\Delta := \overline{M}_{0,d+1}\setminus M_{0,d+1}$ is
formed by the following prime divisors: for each subset $I \subset
\{1,2,\ldots,d+1 \}$ with $|I|\geq 2$ and $|I^c|\geq 2$, we let
$\delta_I \hookrightarrow \overline{M}_{0,d+1}$ be the divisor
whose generic element is a curve with two components: the points
marked by $I$ in one, and the points marked by $I^c$ in the other.
Hence $\delta_I=\delta_{I^c}$, and usually we will assume $d+1 \in
I$ to avoid repetitions. These divisors are smooth, and $\Delta=
\sum \delta_I$ is a simple normal crossing. The variety
$\overline{M}_{0,d+1}$ represents a fine moduli space, proper and
smooth over $\spec(\Z)$. For $i \in \{1,\ldots,d+2 \}$, the $i$-th
forgetful map $\pi_i: \overline{M}_{0,d+2} \rightarrow
\overline{M}_{0,d+1},$ which forgets the $i$-th marked point and
stabilizes, gives a universal family. We will mainly consider
$$\pi_{d+2}: \overline{M}_{0,d+2} \rightarrow
\overline{M}_{0,d+1}.$$ It has $d+1$ distinguished sections
$\delta_{1,d+2}, \ldots, \delta_{d+1,d+2}$, producing the markings
on the parametrized curves.

\begin{defi}
Let $X$ be a nonsingular projective variety, and let $D$ be a
nonsingular divisor in $X$. Let $B$ be a curve in $X$. We say that
$B$ is \textbf{transversal} to $D$ if locally at any $x \in D \cap
B$, the curve $B$ can be factored in $B_1, \ldots, B_n$ distinct
local irreducible curves (branches of $B$ in $\widehat{\O}_{x,X}$)
so that $B_i.D=1$ for every $i$. If $D$ is a sum of nonsingular
divisors $D_j$, we say that $B$ is transversal to $D$ if it is to
each $D_j$. \label{d21}
\end{defi}

Below a well-known property for stable families, coming from the
construction of $\overline{M}_{0,d+1}$.

\begin{lem}
Let $x$ be a $\K$-point in $\overline{M}_{0,d+1}$. Let $B$ be a
local curve passing through $x$, i.e., $B$ is a irreducible curve
defined by functions in $\widehat{\O}_{x,\overline{M}_{0,d+1}}
\simeq \K[[t_1,\ldots,t_{d-2}]]$. Assume $t_1 t_2 \cdots t_k=0$
defines $\Delta= \sum \delta_I$, so $k \leq d-2$, and that
$t_j|_{B}$ is not identically zero for all $1 \leq j \leq k$.
Consider the
commutative diagram $$ \xymatrix{ R \ar[d]_{\rho} \ar[r] & \overline{M}_{0,d+2} \ar[d]_{\pi_{d+2}}   \\
\overline{B} \ar[r]^{\iota} & \overline{M}_{0,d+1} }$$ where
$\iota$ is the composition of the inclusion of $B$ with its
normalization, so $\overline{B}$ is the normalization of $B$, and
$R$ is defined by base change. Then, the surface $R$ is normal,
and can only have singularities of the form $$\spec \
\K[u,v,t]/(uv-t^m)$$ at the nodes of the singular fiber, for some
$m$. Moreover, the surface $R$ is nonsingular if and only if $B$
is transversal to $\Delta$.  \label{l21}
\end{lem}

A brief outline of the proof. Let $X \rightarrow \spec \ \K$ be
the corresponding stable curve over $t_1=t_2=\ldots=t_{d-2}=0$.
Consider the deformation of $X$ as described in
\cite[pp.79--85]{DM}. At a nonsingular point of $X$, we have a
nonsingular point for $R$, so we pay only attention to the nodes
of $X$. Let $y$ be a node of $X$, corresponding to the
intersection of $B$ with $t_i$ for some $i \in \{1, \ldots, k\}$,
i.e., the node $y$ splits $\{1,\ldots,d+1 \}$ in two subsets $I$
and $I^c$, and $t_i=0$ corresponds to $\delta_I$. At the
corresponding point $y$ in $\overline{M}_{0,d+2}$ (over $\K$), the
local rings and the universal map can be written in the projection
form
$$\widehat{\O}_{x,\overline{M}_{0,d+1}} = \K[[t_1,\ldots,t_{d-2}]] \rightarrow \widehat{\O}_{y,\overline{M}_{0,d+2}}=
\K[[u_i,v_i,t_1,t_2,\ldots,t_{d-2}]]/(u_iv_i-t_i),$$ for suitable
variables $u_i,v_i$. Now, the composition of the inclusion of $B$
with its normalization $\iota$ has the form
$\iota(t)=(u_1t^{m_1},u_2t^{m_2},\ldots,u_{d-2}t^{m_{d-2}})$ for
some units $u_i$'s and a local parameter $t$ on $\overline{B}$.
This is because $B$ is not in $t_j=0$ for all $j$. Hence,
$\iota^*(\delta_I)= u_i t^{m_i}$ for $\delta_I=\{t_i=0 \}$. Since
$R$ is defined through the base change by $\iota$, we have the
isomorphism $$\widehat{\O}_{y,R} \simeq \K[[u_i,v_i,t]]/(u_i
v_i-t^{m_i}),$$ and so, $R$ is nonsingular iff $m_i=1$ for all $i
\in \{1,\ldots,k \}$, i.e., transversal to $\Delta$.

The moduli spaces $\overline{M}_{0,d+1}$ have a beautiful
construction, due to Kapranov \cite{Ka1,Ka2}, as iterated blow-ups
of $\P^{d-2}$ (see below). It follows that curves in $M_{0,d+1}$
are strict transforms of curves in $\P^{d-2}$, which are not
contained in a certain fixed hyperplane arrangement $\H_d$. The
following description of these spaces can be found in
\cite{Ka1,Ka2}.

\begin{defi}
A \textbf{Veronese curve} is a rational normal curve of degree
$d-2$ in $\P^{d-2}$, i.e., a curve projectively equivalent to
$\P^1$ in its Veronese embedding. \label{d22}
\end{defi}

It is a classical fact that any $d+1$ points in $\P^{d-2}$ in
general position lie on a unique Veronese curve. The points $P_1,
\ldots,P_{n+2}$ are said to be in general position if no $n+1$ of
them lie in a hyperplane. The main theorem in \cite{Ka1} says that
the set of Veronese curves in $\P^{d-2}$ and its closure are
isomorphic to $M_{0,d}$ and $\overline{M}_{0,d}$ respectively.

\begin{teo}(Kapranov \cite{Ka1})
Take $d$ points $P_1,\ldots,P_d$ of projective space $\P^{d-2}$
which are in general position. Let $V_0(P_1,\ldots,P_d)$ be the
space of all Veronese curves in $\P^{d-2}$ through these $d$
points $P_i$. Consider it as a subvariety in the Hilbert scheme
$\H$ parametrizing all subschemes on $\P^{d-2}$. Then,
\begin{itemize}
\item[1.] We have $V_0(P_1,\ldots,P_d) \cong M_{0,d}$. \item[2.]
If $V(P_1,\ldots,P_d)$ is the closure of $V_0(P_1,\ldots,P_d)$ in
$\H$, then $V(P_1,\ldots,P_d) \cong \overline{M}_{0,d}$. The
subschemes representing limit positions of curves from
$V_0(P_1,\ldots,P_d)$ are, considered together with $P_i$, stable
$d$-pointed curves of genus $0$, which represent the corresponding
points of $\overline{M}_{0,d}$. \item[3.] The analogs of
statements (a) and (b) hold also for Chow variety instead of
Hilbert scheme.
\end{itemize}
\label{t21}
\end{teo}

\begin{teo} (Kapranov \cite{Ka2})
Choose $d$ general points $P_1,\ldots,P_d$ in $\P^{d-2}$. The
variety $\overline{M}_{0,d+1}$ can be obtained from $\P^{d-2}$ by
a series of blow-ups of all the projective spaces spanned by
$P_i$. The order of these blow ups can be taken as follows:
\begin{itemize}
\item[1.] Points $P_1,\ldots,P_{d-1}$ and all the projective
subspaces spanned by them in order of the increasing dimension;
\item[2.] The point $P_d$, all the lines $\langle P_1,P_d \rangle,
\ldots ,\langle P_{d-2},P_d \rangle$ and subspaces spanned by them
in order of the increasing dimension; \item[3.] The line $\langle
P_{d-1},P_d \rangle$, the planes $\langle P_i,P_{d-1},P_d
\rangle$, $i\neq d-2$ and all subspaces spanned by them in order
of the increasing dimension, etc, etc.
\end{itemize}
\label{t22}
\end{teo}

Let us denote the Kapranov's map in Theorem \ref{t22} by
$\psi_{d+1}: \overline{M}_{0,d+1} \rightarrow \P^{d-2}.$

Some conventions and notations for the rest of the paper. Let us
fix $d$ points in general position in $\P^{d-2}$. We take
$P_1=[1:0: \ldots:0]$, $P_2=[0:1:0: \ldots:0]$, ..., $P_{d-1}=[0:
\ldots:0:1]$ and $P_d=[1:1: \ldots:1]$. The symbol $\langle
Q_1,\ldots,Q_r \rangle$ denotes the projective linear space
spanned by the points $Q_i$. Let \begin{center} $\Lambda_{i_1,
\ldots,i_r}= \langle P_j : j\notin \{i_1, \ldots,i_r \} \rangle$
\end{center} where $1\leq r \leq d-1$ and $i_1, \ldots,i_r$ are
distinct numbers, and let $\H_{d}$ be the union of all the
hyperplanes $\Lambda_{i,j}$. Hence, $\Lambda_{i,j}=\{ [x_1:
\ldots: x_{d-1}] \in \P^{d-2}: x_i=x_j \}$ for $i,j\neq d$,
$\Lambda_{i,d}=\{ [x_1: \ldots :x_{d-1}] \in \P^{d-2}: x_i=0 \}$
and \begin{center} $\H_d = \{ [x_1: \ldots: x_{d-1}] \in \P^{d-2}:
\ x_1 x_2\cdots x_{d-1} \prod_{i<j} (x_j-x_i)=0 \}.$ \end{center}

\begin{example}
For $d=4$, Theorem \ref{t22} says that the map $\psi_5:
\overline{M}_{0,5} \rightarrow \P^2$ is the blow-up of $\P^2$ at
the points $P_1=[1:0:0]$, $P_2=[0:1:0]$, $P_3=[0:0:1]$, and
$P_4=[1:1:1]$. The hyperplane arrangement $\H_4$ is given by the
complete quadrilateral
$$x_1x_2x_3(x_1-x_2)(x_1-x_3)(x_2-x_3)=0.$$ The universal map $\pi_{5}: \overline{M}_{0,5} \rightarrow \overline{M}_{0,4}=\P^1$ is induced by the
pencil of conics (Veronese curves in $\P^2$) containing $P_1$,
$P_2$, $P_3$, and $P_4$. \label{e21}
\end{example}

\section{Arrangements of $d$ sections coming from curves in $M_{0,d+1}$}\label{s3}

Let $B$ be an irreducible projective curve in
$\overline{M}_{0,d+1}$ with $B \cap M_{0,d+1}\neq \emptyset$. By
using Kapranov's map $\psi_{d+1}: \overline{M}_{0,d+1} \rightarrow
\P^{d-2}$, this is the same as giving an irreducible projective
curve $A$ in $\P^{d-2}$ not contained in $\H_d$. The proper
transform of $A$ under $\psi_{d+1}$ is $B$.
Consider the base change diagram $$ \xymatrix{ R \ar[d]_{\rho} \ar[r] & \overline{M}_{0,d+2} \ar[d]_{\pi_{d+2}}   \\
\overline{B} \ar[r]^{\iota} & \overline{M}_{0,d+1} }$$ where
$\iota$ is the composition of the inclusion of $B$ with its
normalization. Let us denote $\overline{B}=C$. Notice that the
distinguished sections $\delta_{1,d+1}, \delta_{2,d+1}, \ldots,
\delta_{d+1,d+1}$ of $\pi_{d+2}$ induce $d+1$ sections
$\tilde{S}_{1,d+2},\ldots,\tilde{S}_{d+1,d+2}$ for $\rho$. Also,
by Lemma \ref{l21}, the surface $R$ is a normal projective surface
with only canonical singularities of type $uv=t^m$ for various
integers $m$, and only at nodes of singular fibers. We now resolve
these singularities minimally to obtain a fibration $\tilde{\rho}:
\widetilde{R} \rightarrow C$, so that $\widetilde{R}$ is
nonsingular. Notice that $\tilde{\rho}$ has only reduced trees of
$\P^1$'s as fibers, and it has $d+1$ distinguished sections.

Let $F$ be a singular fiber of $\tilde{\rho}$. Consider the curves
$E$ in $F$ with $E.(F-E)=1$, and which do not intersect the
$(d+1)$-th section (the proper transform of
$\tilde{S}_{d+1,d+2}$). Then, the $E$'s are disjoint with
self-intersection $-1$. We now blow down all of these $E$'s to
obtain a new fibration over $C$ with $d+1$ distinguished sections,
and reduced trees of $\P^1$'s as fibers. If there is a singular
fiber $F$, we repeat the previous procedure. After finitely many
steps, this stops in a fibration $\rho_0: R_0 \rightarrow C$ with
nonsingular fibers, and $d+1$ labeled sections
$\{S_1,S_2,\ldots,S_{d+1}\}$, where $\tilde{S}_{i,d+2}$ is the
proper transform of $S_i$.

\begin{prop}
The fibration $\rho_0: R_0 \rightarrow C$ is isomorphic over $C$
to $\pi: \P_{C}(\L) \rightarrow C$ with $\L \simeq \iota^*
\big(\psi_{d+1}^*(\O_{\P^{d-2}}(1)) \big)$, and so
$\deg(\L)=e=\deg (\psi_{d+1}(B))$. The labeled set
$\{S_1,\ldots,S_d\}$ is a primitive arrangement of $d$ sections.
\label{p31}
\end{prop}

\begin{proof}
By \cite[V, Proposition 2.8]{Ha1}, the ruled surface $\rho_0: R_0
\rightarrow C$ is isomorphic over $C$ to $\P_C(\E) \rightarrow C$,
where $\E$ is a rank two locally free sheaf on $C$ with the
property that $H^0(\E)\neq 0$ but for all invertible sheaves $\M$
on $C$ with $\deg \M < 0$, we have $H^0(\E \otimes \M)=0$. So, we
assume $R_0=\P_C(\E)$. Since, by construction, $S_i.S_{d+1}=0$ for
all $i\neq d+1$, we have that $\E$ is decomposable, say $\E\simeq
\O_C \oplus \L^{-1}$, where $\L$ is unique up to isomorphism.
Moreover, $S_i \sim S_{d+1} + \pi^*(\L)$ for all $i\neq d+1$. In
particular, $0<S_i.S_j=e=\deg \L$ for all $i,j\neq d+1$ and
$S_{d+1}^2=-e$. Notice that $\bigcap_{i=1}^d S_i =\emptyset$.
Hence, $\A=\{S_1,\ldots,S_d \}$ is an arrangement of $d$ sections
of $\pi: \P_{C}(\L) \rightarrow C$.

Observe that $\A$ is primitive because it comes from the
normalization of a curve in $\overline{M}_{0,d+1}$. Assume there
is a morphism of arrangements from $\A$ to $\A'$, with data
$(C',\L',d)$ and map $g: C \rightarrow C'$ (see Definition
\ref{d12}). Then, our map $\iota: C \rightarrow B$ would factor
through $g$, induced by the natural map $\iota': C' \rightarrow
\overline{M}_{0,d+1}$. This is possible only if $\deg g =1$,
because $C$ is the normalization of $B$, and so $\A$ and $\A'$ are
isomorphic arrangements.

Let $c$ be a point in $C$, and consider the fiber $F_c$ of $\pi$.
Let $S_i,S_j \neq S_{d+1}$ be distinct sections which intersect at
a point $P$ in $F_c$. Then, through the description in Lemma
\ref{l21}, it is not hard to see that $$(S_i.S_j)_P=
C_{\text{loc}}. \sum_{\text{all} \ I \ \text{with} \ i,j \in
I\setminus \{d+1\}} \delta_I$$ where $C_{\text{loc}}$ is the
corresponding local branch of $B$ at $\iota(c)$.

Let $\{P_1,\ldots,P_m \}=S_i \cap S_j$. Let $\sigma_j:C
\rightarrow \P_C(\L)$ be the section corresponding to $S_j$. Then,
$\L \simeq \sigma_j^*(\pi^*(\L)\otimes \O_{S_j}) \simeq
\sigma_j^*(\O_{\P_C(\L)}(S_i) \otimes \O_{S_j}) \simeq \O_C
\Big(\sum_{k=1}^m (S_i.S_j)_{P_k} \pi(P_k) \Big),$ and so, because
of the previous formula, $\L \simeq \iota^*\Big( \sum_{\text{all}
\ I \ \text{with} \ i,j \in I\setminus \{d+1\}} \delta_I \Big)$.
Now, by Kapranov's description in Theorem \ref{t22}, we have
\begin{center} $\psi_{d+1}^*(\O_{\P^{d-2}}(1)) \simeq
\sum_{\text{all} \ I \ \text{with} \ i,j \in I\setminus \{d+1\}}
\delta_I, $ \end{center} and so $\L \simeq \iota^*
\big(\psi_{d+1}^*(\O_{\P^{d-2}}(1)) \big)$. This comes from the
pull-back of the hyperplane $\Lambda_{i,j}$. By the projection
formula, we have
\begin{center} $\deg(\L)=B.\psi_{d+1}^*(\O_{\P^{d-2}}(1))={\psi_{d+1}}_*(B).\O_{\P^{d-2}}(1)=\deg{\psi_{d+1}(B)}.$ \end{center}
\end{proof}

\section{The one-to-one correspondences}\label{s4}

Fix an integer $d\geq 3$, and an algebraically closed field $\K$.

\begin{teo} We have $$\bigsqcup_{C,\L} \AA(C,\L,d) \equiv \{
\text{irreducible curves in} \ M_{0,d+1} \},$$ where the disjoint
union is over all nonsingular projective curves $C$ and line
bundles $\L$ on $C$ (both up to isomorphism). This equality gives
a bijection between $\AA(C,\L,d)$ and the set of irreducible
projective curves $B$ in $\overline{M}_{0,d+1}$ with $M_{0,d+1}
\cap B \neq\emptyset$, whose normalization is $C$ and $\L \simeq
\iota^* \big( \psi_{d+1}^*(\O_{\P^{d-2}}(1)) \big)$, where $\iota:
C \rightarrow B$ is the composition of the inclusion of $B$ and
its normalization. \label{t41}
\end{teo}

\begin{proof}
Let $B$ be an irreducible curve in $\overline{M}_{0,d+1}$ with $B
\cap M_{0,d+1} \neq \emptyset$. By Proposition \ref{p31}, $B$
produces a unique element in $\AA \big(\overline{B},\iota^*
(\psi_{d+1}^*(\O_{\P^{d-2}}(1))),d \big)$, where $\iota$ is the
composition of the inclusion of $B$ and its normalization. In this
way, we only need to prove that given $\A \in \AA(C,\L,d)$, there
is an irreducible curve $B$ in $\overline{M}_{0,d+1}$ intersecting
$M_{0,d+1}$ so that $\A$ is induced by $B$ as in Proposition
\ref{p31}.

Let $\A =\{S_1,\ldots,S_d \}$ be a primitive arrangement of $d$
sections of $\pi: \P_C(\L) \rightarrow C$. The section $C_0$ is
denoted by $S_{d+1}$. We repeatedly perform blow-ups at the
intersections of the sections $S_i$ and their proper transforms,
until they are all disjoint. We do this in a minimal way, that is,
given a $(-1)$-curve in a fiber, its blow-down produces an
intersection of the distinguished sections. The corresponding
fibration $\widetilde{T} \rightarrow C$ has $(d+1)$-pointed
semi-stable genus zero curves as fibers. The $d+1$ markings are
produced by intersecting the proper transforms of the sections
$S_i$'s with the fibers. They may fail to be stable exactly
because of the presence of fibers with $\P^1$'s having no
markings, and intersecting the rest of the fiber at two points.
These components form chains of $\P^1$'s, which we blow down to
obtain a $(d+1)$-pointed stable family of genus zero curves $T
\rightarrow C$. Therefore, we have a commutative
diagram $$\xymatrix{ T \ar[d] \ar[r] & \overline{M}_{0,d+2} \ar[d]_{\pi_{d+2}}   \\
C \ar[r] & \overline{M}_{0,d+1} }$$ so that $T \simeq C
\times_{\overline{M}_{0,d+1}} \overline{M}_{0,d+2}$. Notice that
the image of $C$ is a curve $B$ because $\bigcap_{i=1}^d S_i =
\emptyset$ (so not a point), and $B$ intersects $M_{0,d+1}$. Let
$\overline{B}$ be the normalization of $B$, and let $\iota:
\overline{B} \rightarrow \overline{M}_{0,d+1}$ be
the corresponding map. Then, the diagram above factors as $$ \xymatrix{ T \ar[d] \ar[r] & R \ar[d] \ar[r] & \overline{M}_{0,d+2} \ar[d]_{\pi_{d+2}}   \\
C \ar[r]^{f} & \overline{B} \ar[r]^{\iota} & \overline{M}_{0,d+1}
} $$ where $R$ is given by pull-back, and $T \simeq C
\times_{\overline{B}} R$. Let $\widetilde{R}$ be the minimal
resolution of the singularities of $R$. Let us finally consider
the commutative diagram
$$\xymatrix{ T_0 \ar[d] \ar[r] & \widetilde{R} \ar[d] \\ C \ar[r]^f & \overline{B} }$$
where $T_0 \simeq C \times_{\overline{B}} \widetilde{R}$ (it may
be singular). The pull-back of the $d+1$ distinguished sections
are the $d+1$ distinguished sections of $T_0 \rightarrow C$. We
now inductively blow-down all $(-1)$-curves on the fibers of
$\widetilde{R} \rightarrow \overline{B}$ in the following way.

Let $R_i \rightarrow \overline{B}$ be the fibration produced in
the $i$-th step, where $\widetilde{R}=R_0$. Then, $T_i \simeq C
\times_{\overline{B}} R_i$. We obtain the fibration $R_{i+1}
\rightarrow \overline{B}$ through the commutative diagram
below.$$\xymatrix{
& \widetilde{T}_i \ar[d] \ar[dl] & & & \\
T_{i+1} \ar@/_/[rrd]_{\simeq} & T_i = C \times_{\overline{B}} R_i \ar[rrr] \ar[rd] \ar[ddr] & & & R_i \ar[dl] \ar[ddl] \\
& & C \times_{\overline{B}} R_{i+1} \ar[r] \ar[d] & R_{i+1} \ar[d]& \\
& & C \ar[r]^f & \overline{B} &  }$$ Let $E$ be a $(-1)$-curve in
a fiber of $R_i \rightarrow \overline{B}$, and let $P$ be the
point of intersection with the rest of the fiber. Notice that at
least two distinguished sections $U$ and $V$ intersect $E$ (not at
$P$, of course). Let $R_{i+1}$ be the blow-down of $E$, and
$R_{i+1} \rightarrow \overline{B}$ be the corresponding fibration.
Let $Q$ be the pre-image of $P$ and $F$ the pre-image of $E$ in
$T_i$. Notice that $T_i$ may be singular at $Q$, say with a
singularity of type $xy=t^a$. If $a>1$ we resolve $Q$ to get
$\widetilde{T}_i$. Then we define $T_{i+1}$ to be the blow-down of
the total transform of $F$ in $\widetilde{T}_i$ (this is a chain
of $(-1)$-curves). Let us consider $P$ in $R_{i+1}$, and its
pre-image in $C \times_{\overline{B}} R_{i+1}$, say $Q'$. Now, $C
\times_{\overline{B}} R_{i+1}$ is nonsingular at $Q'$, and there
is a morphism $T_{i+1} \rightarrow C \times_{\overline{B}}
R_{i+1}$, which is clearly an isomorphism.

For $R_i$, these procedure is what we have in Section \ref{s3}.
When it stops, say at the $m$-th step, we have $\pi':
R_m=\P_{\overline{B}}(\L') \rightarrow \overline{B}$ with $d+1$
distinguished sections, and $T_m$ nonsingular $\P^1$-bundle over
$C$. Moreover, because of the construction of $T_i$'s, we have
$T_m \simeq \P_C(\L)$, and the arrangement $\A$ is the pull-back
of the one in $\pi': \P_{\overline{B}}(\L') \rightarrow
\overline{B}$. So this is a morphism of arrangements as in
Definition \ref{d12}. But $\A$ primitive implies $\deg f=1$, and
so $C=\overline{B}$.
\end{proof}

For example, one has $\bigsqcup_{C,\L} \AA(C,\L,3) =
\AA(\P^1,\O(1),3)$ which corresponds to the unique curve in
$M_{0,4}=\P^1 \setminus \{[0:1],[1:1],[1:0]\}$. An immediate
corollary is the following.

\begin{cor}
Given a nonsingular projective curve $C$ and a line bundle $\L$ on
$C$, the Kapranov's map $\psi_{d+1}: \overline{M}_{0,d+1}
\rightarrow \P^{d-2}$ gives a one-to-one correspondence between
elements of $\AA(C,\L,d)$ and irreducible projective curves $A$ in
$\P^{d-2}$ not contained in $\H_d$ such that $\overline{A}=C$ and
$\L \simeq \iota^* \big( \psi_{d+1}^* (\O_{\P^{d-2}}(1)) \big)$.
In particular, $\deg A = \deg \L$. \label{c41}
\end{cor}

\begin{cor}
$\AA(\P^1,\O(1),d) \equiv \{ \text{lines in} \ \P^{d-2} \
\text{not in} \ \H_d \}$. \label{c42}
\end{cor}

\begin{proof}
A curve of degree one in $\P^{d-2}$ is a line.
\end{proof}

\begin{cor}
$\AA(\P^1,\O(2),d) \equiv \{ \text{conics in} \ \P^{d-2} \
\text{not in} \ \H_d \}$. \label{c43}
\end{cor}

\begin{proof}
An irreducible curve of degree two is a conic.
\end{proof}

The next two corollaries identify precisely the two distinguished
classes of arrangements in Definition \ref{d14}.

\begin{cor} $$\bigsqcup_{C,\L} \AA_t(C,\L,d) \equiv \{ \text{irreducible projective curves in} \ \overline{M}_{0,d+1} \text{transversal to} \
\Delta\},$$
where the disjoint union is over all nonsingular projective curves
$C$ and line bundles $\L$ on $C$, both up to isomorphism.
\label{c44}
\end{cor}

\begin{proof}
Let $\A \in \AA_t(C,\L,d)$. Let $P$ be a singular point of the
reducible curve defined by $\A$ in $\P_C(\L)$. Let $F_c$ be the
fiber containing $P$. Hence, since $\A$ satisfies (t) in
Definition \ref{d14}, there are two transversal sections $S_i,S_j$
containing $P$, i.e., $(S_i.S_j)_P=1$. Consider the blow-up at
$P$, $\text{Bl}_P(\P_C(\L)) \rightarrow \P_C(\L)$, and let $E$ be
the exceptional curve. Then, $E$ has at least three special
distinct points: the intersections with $\tilde{F}_c$,
$\tilde{S}_i$, and $\tilde{S}_j$ (corresponding proper
transforms). Now, it is clear that the final \textbf{stable}
fibration $R \rightarrow C$ produced from $\A$ has the proper
transform of $E$ as a component of the fiber over $c$. Let
$\tilde{\A}= \{ \tilde{S}_1, \ldots, \tilde{S}_d \}$ be the
collection of proper transforms of $S_i$'s in
$\text{Bl}_P(\P_C(\L))$. Then, $\tilde{\A}$ satisfies property (t)
in Definition \ref{d14} (extending naturally this definition). So,
we repeat the blow-ups until all sections are disjoint (in a
minimal way) to obtain the stable fibration $R \rightarrow C$,
where no blow-downs are needed. Since $R$ is a nonsingular
surface, the curve $\iota(C)$ is transversal to $\Delta$ by Lemma
\ref{l21}.

Now assume $\A$ is not in $\AA_t(C,\L,d)$. Then, there are indices
$i,j$ and a point $P \in S_i \cap S_j$ such that $n=\text{max} \{
(S_i.S_k)_P \ : \ (S_i.S_k)_P \leq (S_i.S_j)_P -1 \} < (S_i.S_j)_P
-1$. Let $(S_i.S_j)_P=m$, so $0\leq n \leq m-2$. We blow up $n$
times the corresponding point in $\tilde{S}_i \cap \tilde{S}_j$
for the successive proper transforms of $S_i$ and $S_j$. Let $X$
be the resulting surface, and $\tilde{P} = \tilde{S}_i \cap
\tilde{S}_j$. Let $E$ be the exceptional curve of the blow-up at
$\tilde{P}$. Then, $E$ has only two special points: the
intersection with the rest of the fiber and with the section
$\tilde{S}_i$. Notice that $\tilde{S}_i \cap \tilde{S}_j \neq
\emptyset$ at $E$. So, in the process to obtain the corresponding
stable fibration $R \rightarrow C$, we need to blow up again at
$\tilde{S}_i \cap \tilde{S}_j$, and so at the end the proper
transform of $E$ will have to be blown down (in order to have a
stable fibration). Therefore, $R$ is singular, and by Lemma
\ref{l21}, $\iota(C)$ is not transversal to $\Delta$. By Theorem
\ref{t41}, we have checked all irreducible curves in $M_{0,d+1}$.
\end{proof}

\begin{cor} $$\bigsqcup_{C,\L} \AA_s(C,\L,d) \equiv \{ \text{irreducible projective curves in}
\ \P^{d-2} \ \text{transversal to} \ \H_d \},$$ where
the disjoint union is over all nonsingular projective curves $C$
and line bundles $\L$ on $C$, both up to isomorphism. \label{c45}
\end{cor}

\begin{proof}
Let $\A \in \AA_s(C,\L,d)$, and consider its stable fibration
$\rho: R \rightarrow C$. We know that the image of $C$ in
$\overline{M}_{0,d+1}$ is transversal to $\Delta$ by the previous
corollary. Let $c \in C$ be a point whose fiber is singular. Then
there exists an element in $\AA(\P^1,\O(1),d)$ that produces the
same fiber. By Corollary \ref{c42}, the set $\AA(\P^1,\O(1),d)$ is
in one-to-one correspondence with lines in $\P^{d-2}$ not in
$\H_d$. Therefore, there is a line in $\P^{d-2}$ not in $\H_d$
passing through $\psi_{d+1}(\iota(c)) \in \P^{d-2}$. This implies
that the image of $\iota(C)$ under $\phi_{d+1}$ is transversal to
$\H_d$. The converse is clear using the same correspondence with
lines.
\end{proof}

\section{Producing explicit primitive arrangements} \label{s5}

In the previous section we classified all arrangements of $d$
sections (and two distinguished subclasses). They are in
one-to-one correspondence with curves in $\P^{d-2}$ outside of a
certain fixed hyperplane arrangement $\H_d$ (Theorem \ref{t41}).
In \cite{U1}, we used this correspondence to explicitly find new
special line arrangements in $\P^2$. For this, we computed the
corresponding line as in Corollary \ref{c42}. In general, it is
hard to present a curve in $\P^{d-2}$ in the form we need to
construct its corresponding arrangement. In this brief section, we
show a simple way to produce arrangements via irreducible curves
in $\P^2$. This is based on \cite[Section 7]{CT1}, where Castravet
and Tevelev describe how to cover $M_{0,d+1}$ with blow-ups of
$\P^2$ at $d+1$ points.

\begin{prop} \cite[Proposition 7.3]{CT1}
Suppose $p_1,\ldots,p_{d+1}$ are distinct points in $\P^2$, and
let $U \subset \P^2$ be the complement to the union of lines
containing at least two of them. The morphism $$\theta: U
\rightarrow M_{0,d+1}$$ obtained by projecting
$p_1,\ldots,p_{d+1}$ from points of $U$ extends to the morphism
$$\theta: Bl_{p_1,\ldots,p_{d+1}}\P^2 \rightarrow \overline{M}_{0,d+1}.$$
If there is no (probably reducible) conic through
$p_1,\ldots,p_{d+1}$ then $\theta$ is a closed embedding. In this
case the boundary divisors $\delta_I$ of $\overline{M}_{0,d+1}$
pull-back as follows: for each line $L_I:= \langle p_i \rangle_{i
\in I} \subset \P^2$, we have $\theta^*(\delta_{I}) =
\widetilde{L}_I$ (the proper transform of $L_I$) and (assuming
$|I| \geq 3$), $\theta^*(\delta_{I\setminus \{k \}})= E_k$, where
$k \in I$ and $E_k$ is the exceptional divisor over $p_k$. Other
boundary divisors pull-back trivially. \label{p51}
\end{prop}

In this way, we have $$ \theta^*(\psi_{d+1}^*( \O(1))) =
(n_{d+1}-1)H - (n_{d+1}-2)E_{d+1} - \sum_{i=1}^d \epsilon_i E_i$$
where $H$ is the class of a general line in $\P^2$, $n_{d+1}$ is
the number of lines in $\P^2$ passing through $p_{d+1}$ and some
other $p_j$, and $\epsilon_i=0$ if there is a $p_k$ in $\langle
p_{d+1}, p_i \rangle$ $k\neq i,d+1$ or $\epsilon_i=1$ otherwise.
Hence, the image of $Bl_{p_1,\ldots,p_{d+1}}\P^2$ under
$\psi_{d+1} \circ \theta$ is a surface $S$ in $\P^{d-2}$ of degree
$2n_{d+1}-3-\sum_{i=1}^d \epsilon_i$, and so $$2 \leq \deg(S) \leq
d-3 .$$ Therefore, $S$ is a surface of minimal degree in some
$\P^{\deg(S)+1} \subset \P^{d-2}$. Thus $S$ is either a rational
normal scroll in $\P^{\deg(S)+1}$ or the Veronese of $\P^2$ in
$\P^5$. Moreover, $S$ is smooth. One can check that $\psi_{d+1}$
blows down certain $d$ $(-1)$-curves in
$Bl_{p_1,\ldots,p_{d+1}}\P^2$ (proper transforms of lines $\langle
p_{d+1},p_i \rangle$ with $\epsilon_i=1$, and $E_i$ with
$\epsilon_i=0$) having as result a Hirzebruch surface $\F_m$,
where $m$ depends on the configuration of points $p_i$ such that
$\epsilon_i=1$.

Given $p_1,\ldots,p_{d+1}$ points in $\P^2$, with no (probably
reducible) conic through them, we consider an irreducible plane
curve $\Gamma$ not included in the union of lines containing
$p_1,\ldots,p_{d+1}$. Then, by Proposition \ref{p51}, we have the
inclusion $\theta: B:=\widetilde{\Gamma} \hookrightarrow
\overline{M}_{0,d+1}$ and so a primitive arrangement $\A$ in
$\AA(\overline{B},\L,d)$ for some line bundle $\L$, by Theorem
\ref{t41}. The line bundle $\L$ depends on the specific
configuration $p_1,\ldots,p_{d+1}$ and the position of the curve
$\Gamma$ with respect to these points. Proposition \ref{p51} gives
a way to precisely see all possible intersections of $\Gamma$ with
$\Delta$, and so this procedure indeed gives an explicit
description of $\A$.

\begin{center}
\begin{figure}[htbp]
\includegraphics[width=12cm]{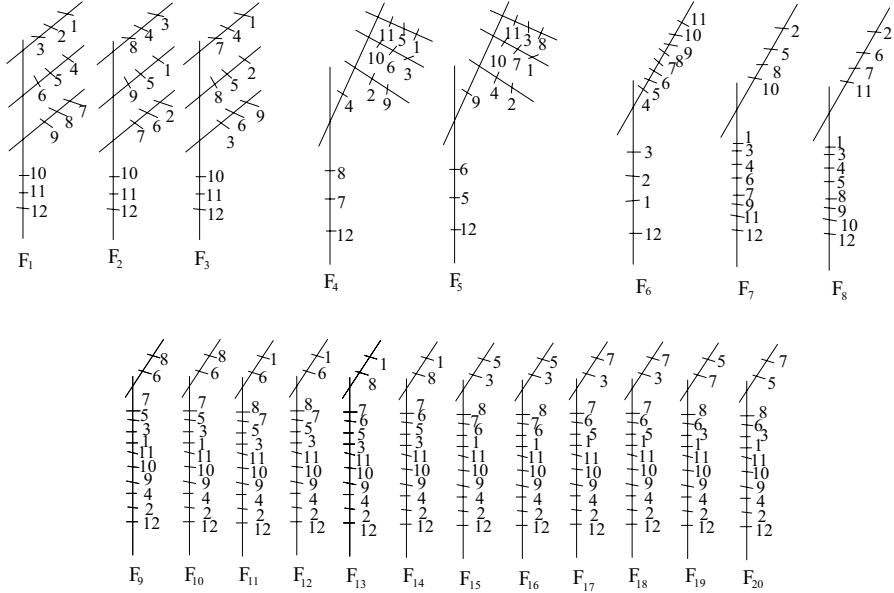} \caption{The singular fibers of the stable fibration induced by the conic $\Gamma$.} \label{f3}
\end{figure}
\end{center}

\begin{example}
Let $\zeta= e^{\frac{2 \pi i}{3}}$. Consider the dual Hesse
arrangement of $9$ lines in $\P^2(\C)$:
$$(x_1^3-x_2^3)(x_1^3-x_3^3)(x_2^3-x_3^3)=0 .$$ It has $12$ triple points and
no other singularities. We label these points as
$p_1=[1:\zeta:\zeta]$, $p_2=[\zeta:\zeta^2:1]$, $p_3=[1:\zeta:1]$,
$p_4=[\zeta:1:\zeta^2]$, $p_5=[\zeta:1:1]$, $p_6=[1:\zeta^2:1]$,
$p_7=[1:1:\zeta]$, $p_8=[1:1:\zeta^2]$, $p_9=[1:1:1]$,
$p_{10}=[0:1:0]$, $p_{11}=[1:0:0]$, and $p_{12}=[0:0:1]$. Consider
the unique conic $\Gamma$ passing through
$p_{12},p_{11},p_{10},p_9,p_4$. It is given by the equation
$\zeta^2 x_1x_2+\zeta x_1x_3+x_2x_3=0$. By Proposition \ref{p51},
we have $\theta: Bl_{p_1,\ldots,p_{12}}\P^2 \hookrightarrow
\overline{M}_{0,12}$, embedding the proper transform
$\widetilde{\Gamma}$ of $\Gamma$. In \cite{CT1}, it is proved that
$\widetilde{\Gamma}$ is a rigid curve in $\overline{M}_{0,12}$.

By Theorem \ref{t41}, this curve defines a primitive arrangement
of $11$ curves $\A$. To actually exhibit $\A$, we need to check
all intersections between $\Gamma$ and all the lines passing
through pairs of points $p_i$ (so, more than the ones in the dual
Hesse arrangement). After that, it is easy to draw a picture of
the arrangement. In Figure \ref{f3}, we show all the singular
fibers of the corresponding stable fibration. Notice that the
arrangement belongs to $\AA_t(\P^1,\O(3),11)$.

For another model, we perform in $\F_3$ two elementary
transformations on the fibers $F_4$ and $F_5$ by blowing up the
corresponding singular points in $\A$. Then, we end up in $\F_1$
where the $(-1)$-curve is the proper transform of $S_{12}$. After
blowing it down, we obtain a very special arrangement of $7$ lines
and $4$ conics in $\P^2$. \label{e51}
\end{example}

\section{Extended and partially extended arrangements} \label{s6}

Fix the data $(C,\L,d)$ over $\K=\overline{\K}$ as in Section
\ref{s1}. We now study properties of certain log surfaces
associated to arrangements of sections $\A$. First, we associate
to $\A$ an extended arrangement $\A_{\Delta}$, and partially
extended arrangements $\A_{p \Delta}$.

\begin{defi}
Consider the arrangement of sections $\A$ as a (reducible
singular) curve, we denote its set of singular points by
$\sing(\A)$. Let $\{ F_1, \ldots, F_{\delta} \}$ be the fibers of
$\pi: \P_C(\L) \rightarrow C$ which contain points in $\sing(\A)$.
Then, the \textbf{extended} arrangement $\A_{\Delta}$ associated
to $\A$ is
$$\A_{\Delta} :=\A \cup \{F_1,\ldots, F_{\delta} \} \cup \{S_{d+1} \}.$$ Let $0 < \epsilon \leq \delta - 2$ be an integer. Let
$\{F_{i_1},\ldots, F_{i_{\epsilon}} \}$ be a subset of $\{ F_1,
\ldots, F_{\delta} \}$ such that for any $1 \leq j \leq \epsilon$
and any point $P$ in $\sing(\A) \cap F_{i_j}$, there are two
sections in $\A$ intersecting at $P$ with distinct tangent
directions. Then, a \textbf{partially extended} arrangement $\A_{p
\Delta}$ associated to $\A$ is
$$\A_{p\Delta} :=\A_{\Delta} \setminus
\{F_{i_1},\ldots,F_{i_{\epsilon}} \}.$$ The numeration of the
fibers will be irrelevant. \label{d61}
\end{defi}

As before, we perform blow-ups at the points in $\sing(\A)$ (and
infinitely near points above them) to separate all sections
$S_i$'s in a minimal way (as in the proof of Theorem \ref{t41}).
This is described by a chain of blow-ups $$ \widetilde{R} = R_s:=
Bl_{P_s} R_{s-1} \stackrel{\rm \sigma_s} {\rightarrow} \cdots
\stackrel{\rm \sigma_3} {\rightarrow} R_{2}:= Bl_{P_2}  R_1
\stackrel{\rm \sigma_2} {\rightarrow} R_1 :=Bl_{P_1} \P_C(\L)
\stackrel{\rm \sigma_1} {\rightarrow} R_0 := \P_C(\L) $$ whose
composition is denoted by $\sigma:\widetilde{R} \rightarrow
\P_C(\L)$. The map $\sigma$ gives the \textbf{minimal log
resolution} of $\A$, and produces the semi-stable fibration of
$(d+1)$-pointed genus zero curves $\tilde{\rho}: \widetilde{R}
\rightarrow C$. The $\sigma^*(\A)_{\text{red}}$ is a simple normal
crossing divisor. Let $t(P_i)$ be the number of sections in the
proper transform of $\A$ passing through $P_i$ right before
blowing up $P_i$ (the center of the blowing up $\sigma_i$). We
define
$$\tau(\A):= \sum_{i=1}^s \big( t(P_i)-1 \big) .$$

The divisor $\sigma^*(\A_{\Delta})_{\text{red}}$ is the minimal
log resolution of $\A_{\Delta}$, but
$\sigma^*(\A_{p\Delta})_{\text{red}}$ may not be minimal, since we
may need to blow down $(-1)$-curves coming from some nodes in
$\A$. The divisors $\sigma^*(\A_{\Delta})_{\text{red}}$ and
$\sigma^*(\A_{p\Delta})_{\text{red}}$ are denoted by
$\bar{\A}_{\Delta}$ and $\bar{\A}_{p\Delta}$ respectively. The
arrangement $\bar{\A}_{\Delta}$ may be seen as defined by the
intersection of the boundary $\Delta$ in $\overline{M}_{0,d+2}$
with the surface $R$, where $\rho: R \rightarrow C$ is the stable
fibration of $(d+1)$-pointed curves of genus zero induced by $\A$.

We now follow the exposition of log surfaces as in \cite[Section
2]{U2}, which is due to Iitaka, and the references given there. We
are interested in the log surfaces
$(\widetilde{R},\bar{\A}_{\Delta})$ and
$(\widetilde{R},\bar{\A}_{p\Delta})$, and their \textbf{log Chern
classes}\footnote{The corresponding log Chern numbers for
$\sigma^*(\A_{p\Delta})_{\text{red}}$ and the minimal log
resolution of $\A_{p\Delta}$ are the same (see for example
\cite[Proposition 2.4]{U2}).}
$$ \bar{c}_i(\bar{\A}_{\Delta}):= c_i
\big(\Omega_{\widetilde{R}}^1(\log \bar{\A}_{\Delta})^{\vee}
\big), \ \ \ \ \ \bar{c}_i(\bar{\A}_{p\Delta}):= c_i
\big(\Omega_{\widetilde{R}}^1(\log \bar{\A}_{p\Delta})^{\vee}
\big)$$ where $\Omega_{\widetilde{R}}^1(\log
\bar{\A}_{\Delta})^{\vee}$, $\Omega_{\widetilde{R}}^1(\log
\bar{\A}_{p\Delta})^{\vee}$ are the dual of the corresponding
sheaves of log differentials (see \cite[Def.2.2]{U2}), and
$i=1,2$.

\begin{prop}
Let $\A$ be an arrangement of sections with data $(C,\L,d)$, $\deg
\L=e$ and $h^1(C,\O_C)=g$. Then, $$\bar{c}_1^2(\A_{\Delta}) =
(d-1)(2\delta + 4(g-1) -e) + \tau(\A) \ \ \ \ \
\bar{c}_2(\A_{\Delta})= (d-1)(2(g-1)+\delta)$$ with $\delta$ as in
Definition \ref{d61}. \label{p61}
\end{prop}

\begin{proof}
In general, if $(Y,D=\sum_{i=1}^r D_i)$ is a log surface (as in
\cite[Section 2]{U2}), the log Chern numbers are (see
\cite[Proposition 2.4]{U2})
$$\bar{c}_1^2(Y,D)= c_1^2(Y) - \sum_{i=1}^r D_i^2 + 2 t_2 + 4 \sum_{i=1}^r (g(D_i)-1),$$ and $\bar{c}_2(Y,D)= c_2(Y) + t_2 + 2 \sum_{i=1}^r
(g(D_i)-1)$, where $t_2$ is the number of nodes of the curve $D$.
In our case, $Y=\widetilde{R}$ and $D=\bar{\A}_{\Delta}$. We will
compute these numbers recursively. Let $\sigma_{i,1}: R_{i}
\rightarrow \P_C(\L)$ be the composition of the blow-ups
$\sigma_{1} \circ \cdots \circ \sigma_{i}$. Define
$$\bar{c}_1^2(i) :=8(1-g) -i - \sum_{j=1}^{s_i} D_{i,j}^2
+2(i+(d+1)\delta) + 4(d+1)(g-1) -4\delta -4i$$ for
$i=0,1,\ldots,s$, where $\sum_{j=1} ^{r_i} D_{i,j}$ is the prime
decomposition of $\sigma_{i,1}^*(\A_{\Delta})_{\text{red}}$. Then,
one can check that $$\bar{c}_1^2(0)=(d-1)(2\delta + 4(g-1)-e) \ \
\ \ \ \ \ \ \bar{c}_1^2(s)= \bar{c}_1^2(\A_{\Delta})$$ and
$\bar{c}_1^2(i+1)= \bar{c}_1^2(i)+ t(P_{i+1})-1$. Therefore,
$\bar{c}_1^2(\A_{\Delta}) = (d-1)(2\delta + 4(g-1) -e) +
\tau(\A)$. On the other hand, by the formula for
$\bar{c}_2(\tilde{R},\bar{\A}_{\Delta})$ above, we have
$$\bar{c}_2(\A_{\Delta})= 4(1-g) + s + (s+(d+1)\delta) +2(d+1)(g-1) - 2\delta - 2s =(d-1)(2(g-1)+\delta). $$
\end{proof}

\begin{cor}
Let $\A=\{S_1,\ldots,S_d \}$ as above. Then,
$\bar{c}_1^2(\A_{\Delta})\geq 2d-1$, $\bar{c}_2(\A_{\Delta})\geq
d-1$, and $2  <
\frac{\bar{c}_1^2(\A_{\Delta})}{\bar{c}_2(\A_{\Delta})}$. The log
canonical divisor $K_{\tilde{R}} + \bar{\A}_{\Delta}$ is big and
nef, and so the surface $\tilde{R}\setminus \bar{\A}_{\Delta}$ is
of log general type. When $\K=\C$, we have the (strict) log
Miyaoka-Yau inequality
$$\bar{c}_1^2(\A_{\Delta}) < 3 \bar{c}_2(\A_{\Delta}),$$ and so $\tau(\A) < (d-1)(\delta + 2(g-1)+e)$. \label{c61}
\end{cor}

\begin{proof}
The second log Chern number
$\bar{c}_2(\A_{\Delta})=(d-1)(2(g-1)+\delta)$ is positive because
$\delta \geq 3$. For the other one, take $S_1 \in \A$. Then, by
looking at the intersections of $S_1$ with $S_2,\ldots,S_d$ in
$\P_C(\L)$, one sees that $\tau(\A)>e(d-1)$. The inequality is
strict because $\bigcap_{i=1}^d S_i= \emptyset$. Therefore, by
Proposition \ref{p61}, one has $\bar{c}_1^2(\A_{\Delta})>0$. By
the same reason,
$$2  < \frac{\bar{c}_1^2(\A_{\Delta})}{\bar{c}_2(\A_{\Delta})} = 2 + \frac{\tau(\A)-e(d-1)}{(d-1)(2(g-1)+\delta)}.$$

Let $D:= \bar{\A}_{\Delta}$ and write its prime decomposition $D=
\sum_{i=1}^{d+1} \tilde{S}_i + \sum_{i} E_i$ where $\tilde{S}_i$
is the proper transform of $S_i$ under $\sigma: \tilde{R}
\rightarrow \P_C(\L)$, and $E_i \simeq \P^1$'s are the rest. Since
$\tilde{S}_{d+1}=\sigma^* S_{d+1}$, we denote $S_{d+1}=\sigma^*
S_{d+1}$ and $F= \sigma^* F$ where $F$ is a general fiber of $\pi:
\P_C(\L) \rightarrow C$. Then, $K_{\tilde{R}} + D \equiv -2
S_{d+1} + (2g-2-e)F + \sum_{i} a_i E_i + D$ for some $a_i >0$.
But, for any fixed $j=1,\ldots,d$, $\tilde{S}_j \equiv S_{d+1} +
eF - \sum_{i} b_i E_i$ where $b_i \leq a_i$, and $\sum_{i} E_i
\equiv \delta F$. Thus, say for $j=1$, we have
$$K_{\tilde{R}} + D \equiv (2g-2+\delta)F + \sum_{i} (a_i-b_i) E_i + \sum_{i=2}^{d} \tilde{S}_i,$$ which means that the log canonical class is
numerically equivalent to an effective divisor. Moreover,
$\tilde{S}_i.(K_{\tilde{R}} + D) = 2g-2 + \delta >0$ for all $i$,
and $E_i.(K_{\tilde{R}} + D)=-2+ E_i.(D-E_i) \geq 0$, and so
$K_{\tilde{R}} + D$ is nef. This plus the fact
$\bar{c}_1^2(\A_{\Delta})=(K_{\tilde{R}}+D)^2>0$ implies that
$K_{\tilde{R}}+D$ is big, and so $\tilde{R} \setminus
\bar{\A}_{\Delta}$ is of log general type. When $\K=\C$,  we use
Sakai's theorem \cite[Theorem 7.6]{Sakai1} to conclude
$\bar{c}_1^2(\A_{\Delta}) \leq 3 \bar{c}_2(\A_{\Delta})$. Notice
that the curve $D$ is semi-stable and has no exceptional curves
with respect to $D$, as defined in \cite[pp.90-91]{Sakai1}. For
strictness, we apply Lemma \ref{apl}.
\end{proof}

\begin{obs}
From the previous proof, it is easy to see that the log canonical
class is ample if and only if $\A$ is transversal (Definition
\ref{d14}). One just applies the Nakai-Moishezon criterion
\cite[p.365]{Ha1}. \label{o61}
\end{obs}

\begin{prop}
Let $\A$ be an arrangement of sections with data $(C,\L,d)$, and
let $\A_{p\Delta}$ be a partially extended arrangement with
$\{F_{i_1},\ldots,F_{i_{\epsilon}} \}$ for some $0<\epsilon \leq
\delta-2$ (see Definition \ref{d61}). Let $k_j^o:= |\sing(\A) \cap
F_{i_j}|$ and $k_j:= |\A \cap F_{i_j}|+1 \leq d$. Then,
$$\bar{c}_2(\A_{p\Delta}) = \bar{c}_2(\A_{\Delta}) - \sum_{j=1}^{\epsilon} k_j + 2 \epsilon \ \ \
\ \ \bar{c}_1^2(\A_{p\Delta})= \bar{c}_1^2(\A_{\Delta}) -
\sum_{j=1}^{\epsilon} k_j^o - 2 \sum_{j=1}^{\epsilon} k_j + 4
\epsilon.$$ \label{p62}
\end{prop}

\begin{proof}
The result follows directly from the formulas in \cite[Proposition
2.4]{U2}.
\end{proof}

\begin{cor}
Let $\A$ be as above. Then, $\bar{c}_1^2(\A_{p\Delta})\geq 2$,
$\bar{c}_2(\A_{p\Delta})\geq 1$. The corresponding log canonical
divisor is big and nef, and so the log surface defined by
$\bar{\A}_{p\Delta}$ is of log general type. When $\K=\C$, we
again have (strict) log Miyaoka-Yau inequality
$\bar{c}_1^2(\A_{p\Delta}) < 3 \bar{c}_2(\A_{p\Delta})$.
\label{c62}
\end{cor}

\begin{proof}
Notice that $\sum_{j=1}^{\epsilon} k_j - 2 \epsilon \leq d
\epsilon - 2 \epsilon \leq (d-2)(\delta-2)$, and so $$
\bar{c}_2(\A_{p\Delta}) = \bar{c}_2(\A_{\Delta}) - \big(
\sum_{j=1}^{\epsilon} k_j - 2 \epsilon \big) \geq
(d-1)(2(g-1)+\delta) - (d-2)(\delta-2) = 2g(d-1) + \delta -2 \geq
1.$$

Clearly, we have $k_i^o + 2 k_i \leq 2d+1$, and so
$\sum_{j=1}^{\epsilon} k_j^o + 2 \sum_{j=1}^{\epsilon} k_j - 4
\epsilon \leq (2d-3)(\delta-2)$. Then, since $\tau(\A)-e(d-1) \geq
1 $, we have $$\bar{c}_1^2(\A_{p\Delta})= \bar{c}_1^2(\A_{\Delta})
- \sum_{j=1}^{\epsilon} k_j^o - 2 \sum_{j=1}^{\epsilon} k_j + 4
\epsilon \ \ \ \ \ \ \ \ \ \ \ \ \ \ \ \ \ \ \ \ \ \ \ \ \ \ \ \ \
\ \ \ \ \ \ \ \ \ \ \ \ \ \ \ \ \ \ \ \ \ \ \ \ \ \ \ \ \ \ \ \ \
\ \ $$
$$ \geq 4(g-1)(d-1)+2\delta(d-1) -e(d-1) + \tau(\A) - (2d-3)(\delta-2) \ \ \ \ \ \ \ \ \ \ $$
$$ \ \ \ \ \ \ \ \ \ \ \ \  \geq 1 + 4(g-1)(d-1) + 2(2d-3) + \delta \geq 1- 4(d-1) + \delta + 2(2d-3) = \delta -1 \geq 2 .$$

We prove nefness and bigness as we did in Corollary \ref{c61}. It
is enough to do it in $\tilde{R}$. Let $D:=
\sigma^*(\A_{p\Delta})_{\text{red}}$ and write its prime
decomposition $D= \sum_{i=1}^{d+1} \tilde{S}_i + \sum_{i} E_i$
where $\tilde{S}_i$ is the proper transform of $S_i$ under
$\sigma: \tilde{R} \rightarrow \P_C(\L)$, and $E_i \simeq \P^1$'s
are the rest. Since $\tilde{S}_{d+1}=\sigma^* S_{d+1}$, we denote
$S_{d+1}=\sigma^* S_{d+1}$ and $F= \sigma^* F$ where $F$ is a
general fiber of $\pi: \P_C(\L) \rightarrow C$. Then,
$K_{\tilde{R}} + D \equiv -2 S_{d+1} + (2g-2-e)F + \sum_{i} a_i
E_i + D$ for some $a_i >0$. But, for any fixed $j=1,\ldots,d$,
$\tilde{S}_j \equiv S_{d+1} + eF - \sum_{i} b_i E_i$ where $b_i
\leq a_i $, and $\sum_{i} E_i \equiv \delta F -
\sum_{j=1}^{\epsilon} \tilde{F}_{i_j}$ (see Definition \ref{d61}).
Thus, say for $j=1$, we have $$K_{\tilde{R}} + D \equiv
(2g-2+\delta - \epsilon)F + \sum_{i} (a_i-b_i) E_i +
\sum_{i=2}^{d} \tilde{S}_i + \sum_{j=1}^{\epsilon}
(F-\tilde{F}_{i_j}) ,$$ which means that the log canonical class
is numerically equivalent to an effective divisor. Moreover,
$\tilde{S}_i.(K_{\tilde{R}} + D) \geq 2g-2+ \delta - \epsilon \geq
0$ for all $i$, and $E_i.(K_{\tilde{R}} + D)=-2+ E_i.(D-E_i) \geq
0$ (this is by definition of partially extended arrangement), and
so $K_{\tilde{R}} + D$ is nef. This plus the fact
$(K_{\tilde{R}}+D)^2>0$ implies that $K_{\tilde{R}}+D$ is big, and
so $\tilde{R} \setminus \bar{\A}_{p\Delta}$ is of log general
type. When $\K=\C$ and by Sakai's theorem \cite[Theorem
7.6]{Sakai1}, we have $\bar{c}_1^2(\A_{p\Delta}) \leq 3
\bar{c}_2(\A_{p\Delta})$. For strictness, we apply again Lemma
\ref{apl}.
\end{proof}

\begin{example} We use Example \ref{e51} to show differences in $\frac{\bar{c}_1^2}{\bar{c}_2}$
between extended and partially extended arrangements. Let $\A$ be
the arrangement in Example \ref{e51}. For the partially extended
arrangements, let $\Xi$ be the set of fibers we take out from
$\A_{\Delta}$ (see Definition \ref{d61}). We label fibers
according to Figure \ref{f3}. Then we have the following table.
\bigskip

\begin{center}
\begin{tabular}{|c|c|c|c|c|c|c|}
\hline

$\Xi$ & $\{F_1,\ldots,F_8 \}$ & $\emptyset$ & $\{F_9,\ldots,F_{20}
\}$ & $\{F_7,\ldots,F_{20}\}$ & $\{F_4,\ldots,F_{20}\}$ &
$\{F_6,\ldots,F_{20}\}$  \\ \hline
 $\bar{c}_1^2$ & $319$ & $399$ & $171$ & $141$ & $124$ & $134$ \\ \hline
 $\bar{c}_2$ & $147$ & $180$ & $72$ & $58$ & $51$ & $55$ \\ \hline
{\tiny $\frac{\bar{c}_1^2}{\bar{c}_2}$} & $2.170...$ &
$2.21\overline{6}$ & $2.375$ & $2.4310...$ & $2.4313...$ &
$2.4\overline{36}$ \\ \hline
\end{tabular}
\end{center}

\label{e61}
\end{example}

\begin{obs}
In general, there are no inequalities between
$\frac{\bar{c}_1^2(\A_{p\Delta})}{\bar{c}_2(\A_{p\Delta})}$ and
$\frac{\bar{c}_1^2(\A_{\Delta})}{\bar{c}_2(\A_{\Delta})}$. Also,
in general, we do not have $2<
\frac{\bar{c}_1^2(\A_{p\Delta})}{\bar{c}_2(\A_{p\Delta})}$. For
instance take a general arrangement $\A \in \AA(\P^1,\O(1),d)$ (a
general line arrangement). We have $\delta=  {d \choose 2}$. Take
$\epsilon=\delta-2$. Then, $\tau(\A)= {d \choose 2}$ and
$\frac{\bar{c}_1^2(\A_{p\Delta})}{\bar{c}_2(\A_{p\Delta})}= 2 -
\frac{d-3}{{d \choose 2} -2}<2$. \label{o63}
\end{obs}

\begin{obs}
Almost any line arrangement is ``log equivalent" to a
$\A_{p\Delta} \in \AA(\P^1,\O(1),d)$. More precisely, let $\LL$ be
a line arrangement in $\P^2$. Assume this line arrangement has at
least two singular points so that each of them belongs to more
than two lines (in particular, general line arrangements are not
allowed). Take two distinct lines in $\LL$, each of which contains
exactly one of these two points. Now we blow up the intersection
of these two lines. The total (reduced) transform of $\LL$ is a
$\A_{p\Delta}$ for some $\A \in \AA(\P^1,\O(1),d)$.

Over $\C$ and by \cite[Theorem 7.2]{U2}, we have the
Hirzebruch-Sakai inequality
\begin{center} $\bar{c}_1^2(\A_{p\Delta}) \leq \frac{8}{3} \,
\bar{c}_2(\A_{p\Delta}), $
\end{center} for any $\A \in \AA(\P^1,\O(1),d)$, with equality if
and only if $\A_{p\Delta}$ is the proper transform of the dual
Hesse arrangement (see Example \ref{e51}). Of course, the same
holds true for $\A_{\Delta}$, but now equality is impossible. This
shows that the log Miyaoka-Yau inequality in the previous
corollaries may be far from optimal, and opens the interesting
question:

\begin{question}
Find the number $0<\alpha(C,\L)< 3$ which makes
$${\bar{c}_1^2(\A_{p\Delta})} \leq \alpha(C,\L) \, {\bar{c}_2(\A_{p\Delta})}$$ a sharp inequality, valid for any complex arrangement $\A \in
\AA(C,\L,d)$ and any $d\geq3$.
\end{question}

Arrangements holding equality should be interesting. Proposition
\ref{p63} will show that we indeed need to look only at
arrangements in $\AA(C,\L,d)$, i.e. primitive ones (Definition
\ref{d13}). \label{o64}
\end{obs}

\begin{obs} The log Miyaoka-Yau inequalities in Corollaries \ref{c61}
and \ref{c62} are not combinatorial, except in the case of
$\AA(\P^1,\O(1),d)$. As in the previous remark, let $\A_{p\Delta}$
be defined by an arrangement of lines $\LL=\{L_1,\ldots,L_s\}$ in
$\P^2$. Let $\{p_1,\dots,p_r\}$ be the set of singular points of
$\LL$. The log Miyaoka-Yau inequality $\bar{c}_1^2(\LL) \leq 3
\bar{c}_2(\LL)$ becomes precisely $r \geq s$ \cite{U2}. This
inequality was proved in a purely combinatorial manner by N. G. de
Bruijn and P. Erd$\ddot{o}$s in \cite{BruErd1948}. Moreover, they
show that $s=r$ if and only if $\LL$ has either $s-1$ lines
through a common point (in this case $\bar{c}_1^2=\bar{c}_2=0$) or
it is a finite projective plane. We proved this inequality in
\cite[Theorem 7.2]{U2} using some surface theory. The following is
another combinatorial proof. Consider $L_i$ as a vector in $\Q^r$
having a $1$ in the $j$-th coordinate if $p_j \in L_i$, $0$
otherwise. The assertion follows if $\LL$ forms a linearly
independent set. If not, say the line $L_1=\sum_{i=2}^s a_i L_i$.
Then, consider the inner product with $L_1-L_j$, giving
$a_j=\frac{L_1\cdot L_1 -1}{1-L_j\cdot L_j} <0$ for all $j$, a
contradiction.

In general, one can exhibit ``combinatorial arrangements" for
which the inequality $\tau(\A) \leq (d-1)(\delta + 2(g-1)+e)$ does
not hold. For example, one may take $\A \in \AA(\P^1,\O(e),4)$
with $e>1$ such that any two sections are tangent of order $e$ and
$\delta=3$ (when $e=1$, $\A_{\Delta}$ is the Fano arrangement
(with a point blown-up)). This combinatorial phenomena is produced
by the freedom we have with respect to tangencies of higher order.
See also Example \ref{o65}, where positive characteristic is used.

\label{o62}
\end{obs}

\begin{prop}
Let $\A,\A'$ be two arrangements of $d$ sections so that $\A$ is a
pull-back of $\A'$, as in Definition \ref{d12}. Assume that the
pull-back map $g$ is separable. Then,
$$\frac{\bar{c}_1^2(\A_{\Delta})}{\bar{c}_2(\A_{\Delta})} \leq
\frac{\bar{c}_1^2(\A'_{\Delta})}{\bar{c}_2(\A'_{\Delta})}, \ \ \ \
\ \ \frac{\bar{c}_1^2(\A_{p\Delta})}{\bar{c}_2(\A_{p\Delta})} \leq
\frac{\bar{c}_1^2(\A'_{p\Delta})}{\bar{c}_2(\A'_{p\Delta})},$$
where $\A_{p\Delta}$ is the pull-back of $\A'_{p\Delta}$.
\label{p63}
\end{prop}

\begin{proof}
By definition, we have the following working diagram
\begin{center} $ \xymatrix{ \P_C(\L) \ar[d]_{\pi} \ar[r]^{G} & \P_{C'}(\L') \ar[d]_{\pi'}   \\ C
\ar[r]^{g} & C' }$ \end{center} where $\A$ and $\A'$ are
arrangements of $d$ sections in $\P_C(\L)$ and $\P_{C'}(\L')$
respectively. Since $g$ is a finite separable morphism, we have
the Hurwitz formula $$2g-2=\deg(g)(2g'-2) + \deg R$$ where
$R=\sum_{c \in C} \text{length}(\Omega_{C/C'})_P P$, and so $\deg
R \geq \sum_{c \in C} (e_c -1)$ (see \cite[p.301]{Ha1}). Here
$e_c$ is the ramification index of $g$ at $c$. As usual, $c' \in
C'$ is a branch point of $g$ if there is $c \in C$ such that
$g(c)=c'$ and $e_c>1$. We remak that for any $c' \in C'$ we have
$g^*(c')= \sum_{g(c)=c'} e_c c$ and $\deg(g)=\sum_{g(c)=c'} e_c$
\cite[p.138]{Ha1}.

It is not hard to see that $\delta + \deg R = \deg(g) \delta' +
\alpha$, for some integer $\alpha \geq 0$. We also have $e=\deg(g)
e'$ (Lemma \ref{l11}, where $e= \deg \L$ and $e'=\deg \L'$).
Notice that, by definition, the map $g$ cannot be branched at any
of the images of the special $\epsilon$ fibers (this is an empty
statement when we consider the extended arrangement). This is
because a pre-image of such a fiber would contain at least one
point in $\sing(\A)$ where all sections in $\A$ through it have
the same tangent direction. So, $\epsilon= \deg(g) \epsilon'$,
$\sum_{i=1}^{\epsilon} k_i= \deg(g) \sum_{i=1}^{\epsilon'} k'_i$,
and $\sum_{i=1}^{\epsilon} k_i^o = \deg(g) \sum_{i=1}^{\epsilon'}
{k_i^o}'$. So we only need to compare $\tau(\A)$ with $\tau(\A')$.

For any $P \in \sing(\A)$, define $t_P(\A):=\sum_{Q \in N(P)}
(t(Q)-1)$ where $N(P)$ is the set of points blown up by $\sigma$
above $P$ (so $N(P)$ contains $P$). Then, $t_P(\A)= e_{\pi(P)}
t_{G(P)}(\A')$, and so $\tau(\A) = \deg(g) \tau(\A')$. Therefore,
$$ \ \ \ \ \ \ \ \frac{\bar{c}_1^2(\A_{p\Delta})}{\bar{c}_2(\A_{p\Delta})} = 2 + \frac{\tau(\A)-e(d-1)-\sum_{i=1}^{\epsilon}k_i^o}{(d-1)(\delta+2(g-1)) - \sum_{j=1}^{\epsilon} k_j + 2 \epsilon} \ \ \ \ \ \ \ \ \ \ \ \ \ \ \ \ \ \ \ \ \ \ \ \ \ \ \ \ \ \ \ \ \ \ \ \ \ \ \ \ \ \ \ \ \ \ \ \  $$
$$  = 2 + \frac{\deg(g)(\tau(\A')-e'(d-1)- \sum_{i=1}^{\epsilon}{k_i^o}') }{(d-1)( \deg(g)(2g'-2 + \delta') + \alpha)
- \deg(g)(\sum_{j=1}^{\epsilon'} k'_j - 2 \epsilon') }$$ $$ \ \ \
\ \ \ \ \ \   \leq 2 + \frac{\deg(g)(\tau(\A')-e'(d-1)-
\sum_{i=1}^{\epsilon}{k_i^o}') }{(d-1) \deg(g)(2g'-2 + \delta')  -
\deg(g)(\sum_{j=1}^{\epsilon'} k'_j - 2 \epsilon') } =
\frac{\bar{c}_1^2(\A'_{p\Delta})}{\bar{c}_2(\A'_{p\Delta})}.$$
\end{proof}

\begin{obs}
The situation is different when the base change is not separable.
Assume that $\K$ has positive characteristic $p$. Let $\A,\A'$ be
two arrangements of $d$ sections so that $\A$ is a pull-back of
$\A'$, as in Proposition \ref{p63}, but now let $g: C=C'_{p^r}
\rightarrow C'$ be the composition of the $\K$-linear Frobenius
morphism $r$ times \cite[p.302]{Ha1}. Then,
$$\frac{\bar{c}_1^2(\A_{\Delta})}{\bar{c}_2(\A_{\Delta})} = 2 + p^r \frac{t(\A')-e'(d-1)}{(d-1)(2(g'-1)+\delta')} = 2 + p^r
\Big(\frac{\bar{c}_1^2(\A'_{\Delta})}{\bar{c}_2(\A'_{\Delta})} -2
\Big)$$ and so it becomes arbitrarily large when $r>>0$ (Corollary
\ref{c61}). In the next section, these examples will produce
nonsingular projective surfaces violating any Miyaoka-Yau
inequality. \label{o65}
\end{obs}

\section{Random surfaces associated to arrangements} \label{s7}

Fix the data $(C,\L,d)$ over $\K=\overline{\K}$ as always. We now
associate to each arrangement of sections $\A$ of $\pi: \P_{C}(\L)
\rightarrow C$ various collections of nonsingular projective
surfaces. Each collection is produced by either $\A_{\Delta}$ or
some $\A_{p\Delta}$. The construction is analogue to the one in
\cite[Theorem 6.1]{U2} but now we have more singular arrangements
of curves. From now on, we consider $\A_{\Delta}$ as
$\A_{p\Delta}$ with $\epsilon=0$, to save notation.

Let $\A=\{S_1,\ldots,S_d \}$. By definition, $$S_i \sim S_{d+1} +
\pi^*(\L)$$ for $i=1,\ldots,d$, where $S_{d+1}=C_0$ with
$C_0^2=-e=-\deg \L$. Let $\{F_1,\ldots,F_{\delta-\epsilon} \}$ be
the fibers which define $\A_{p\Delta}=\A \cup
\{F_1,\ldots,F_{\delta -\epsilon},S_{d+1} \}$. Let $\{0<x_i
\}_{i=1}^d$, $\{0<y_i \}_{i=1}^{\delta-\epsilon}$ be an integer
solution of the equation $$\EE: \ \ \sum_{i=1}^d e x_i +
\sum_{i=1}^{\delta -\epsilon} y_i =p
$$ for some prime number $p\neq \text{char}( \K )$, and let $x_{d+1}:=p - \sum_{i=1}^d x_i$. When $p$ is large enough, the equation $\EE$ has
nonnegative solutions, exactly (see \cite{BGK01}) $$
\frac{p^{d+\delta-\epsilon-1}}{(d+\delta-\epsilon-1)! e^d} +
O(p^{d+\delta-\epsilon-2}). $$ In this way,
$$\sum_{i=1}^{d+1} x_i S_i + \sum_{i=1}^{\delta-\epsilon} y_i F_i \sim p S_{d+1} + \pi^* \Big( \big(\sum_{i=1}^d x_i \big) \L + \M \Big),  $$
for some line bundle $\M$ on $C$ of degree
$\sum_{i=1}^{\delta-\epsilon} y_i$. Then, since $\big(\sum_{i=1}^d
x_i \big)e + \deg \M =p$, there is an invertible sheaf $\N$ on $C$
such that $$\sum_{i=1}^{d+1} x_i S_i +
\sum_{i=1}^{\delta-\epsilon} y_i F_i \sim p ( S_{d+1} + \pi^*
\N).$$

The theorem below associates to each arrangement $\A$ various
families of random smooth projective surfaces. We use the method
in \cite{U2}, with an extra care of the new singularities; in
\cite{U2} we only had simple crossings (as in Definition
\ref{d14}). The randomness part relies on a large scale behavior
of Dedekind sums and continued fractions (see \cite[Apendix]{U2}).
The proof will be based on the work done in \cite{U2}.

\begin{teo} Let $\A$ be an arrangement of sections of $\pi:\P_C(\L) \rightarrow C$. Then, there exist nonsingular projective surfaces $X$ of
general type with

\vspace{0.2 cm}
\begin{center} {\Large$\frac{c_1^2(X) }{c_2(X)}$} arbitrarily close to {\Large $\frac{\bar{c}_1^2(\A_{p\Delta})}{\bar{c}_2(\A_{p\Delta})}$}, \end{center}
for any $\A_{p\Delta}$. \label{t71}
\end{teo}

\begin{proof}
Let $Z:=\P_C(\L)$ and let $Y$ be the surface which log resolves
minimally the arrangement $\A_{p\Delta}$ (for example,
$Y=\tilde{R}$ when $\epsilon=0$). Let $\sigma: Y \rightarrow Z$ be
the minimal log resolution of $\A_{p\Delta}$. Choose a solution of
$\EE$, and define $$D:= \sigma^* \Big(\sum_{i=1}^{d+1} x_i S_i +
\sum_{i=1}^{\delta-\epsilon} y_i F_i \Big) \sim p \ \sigma^*(
S_{d+1} + \pi^* \N). $$ This allows the construction of the $p$-th
root cover $f: X \rightarrow Y$ along $D$, as in \cite[Section
2]{U2}. Thus, $X$ is a nonsingular projective surface. Let $$D=
\sum_{i=1}^r \nu_i D_i$$ be the decomposition of $D$ into the sum
of prime divisors. From $\EE$ and the nature of $\A_{p\Delta}$,
ones sees that $0< \nu_i <p$.

As in \cite[Apendix]{U2}, for $0<q<p$, we denote the corresponding
Dedekind sum by $s(q,p)$ and the length of the corresponding
negative continued fraction by $l(q,p)$. In \cite[Proposition 3.4,
3.6, and 2.4]{U2}, we computed the Chern numbers of $X$ as
functions of $p$, Chern numbers, and ``error terms". Let
$\bar{c}_1^2:= \bar{c}_1^2(\A_{p\Delta})$,
$\bar{c}_2:=\bar{c}_2(\A_{p\Delta})$, $c_1^2:= c_1^2(Y)$, and
$c_2:= c_2(Y)$. Then,
$$c_1^2(X)= \bar{c}_1^2 p + 2(c_2-\bar{c}_2) + (c_1^2-\bar{c}_1^2+ 2\bar{c}_2- 2 c_2)\frac{1}{p} - \sum_{i<j} c(p-\nu_i' \nu_j,p)  D_i.D_j  $$
$$c_2(X)= \bar{c}_2 p + (c_2-\bar{c}_2) + \sum_{i<j} l(p-\nu_i' \nu_j,p) D_i.D_j $$ where $c(p-\nu_i' \nu_j,p):=12s(p-\nu'_i \nu_j,p) + l(p-\nu'_i
\nu_j,p)$. Let us denote the error terms by $CCF := \sum_{i<j}
c(p-\nu_i' \nu_j,p) D_i.D_j$ and $LCF:= \sum_{i<j} l(p-\nu_i'
\nu_j,p) D_i.D_j$.

We prove the existence of ``good" solutions of $\EE$ for
arbitrarily large primes $p$, which make $\frac{CCF}{p}$ and
$\frac{LCF}{p}$ arbitrarily small. In addition, this will show
that random partitions are ``good", with probability approaching
$1$ as $p$ becomes arbitrarily large. The key numbers to take care
of are the $p - \nu_i' \nu_j$, which are defined for each node of
$D_{\text{red}}$. The idea is to show that there are solutions of
$\EE$ for which all $p - \nu_i' \nu_j$ are outside of a certain
\textbf{bad set} $\FF \subset \{0, \ldots , p-1 \}$ (defined in
\cite[Apendix]{U2}, due to K. Girstmair).

We write down for each node in $D_{\text{red}}$ the multiplicities
$\nu_a,\nu_b$ as functions on the numbers $x_i,y_j$. There are
different cases, all described in the following table.

\begin{center}
\begin{tabular}{|c|c|c|c|c|c|}
\hline

\text{Type} & I & II & III & IV & V  \\ \hline
 $\nu_a$ & $x_i$ & $y_i$ & $y_i$ & $\sum_{k=1}^d n_k x_k + y_j$  & $n \sum_k x_k + z,$ $z\neq 0$  \\
  & & & & with $0\leq n_k \leq e$ & has no $x_k$, $0\leq n <e$ \\ \hline
 $\nu_b$ & $x_j$ & $x_j$ & $x_{d+1}$ & $ x_k$ with $n_k \neq 0$  & $\sum_k x_k + \nu_a$  \\ \hline

\end{tabular}
\end{center}

Notice that ``$z\neq 0$ has no $x_k$" in case V because of our
restriction on tangent directions at the singular points of
$\A_{p\Delta}$ (Definition \ref{d61}). Below we estimate for each
type the number of solutions $b(\nu_a,\nu_b)$ of $\EE$ producing a
bad multiplicity $p-\nu_a'\nu_b \in \FF$. We do it case by case.

(\textbf{Type I)} This is a node in $S_i \cap S_j$ (possible only
when $\epsilon>0$). Since $\EE$ is a weighted partition of $p$, we
can use the estimate in \cite[proof of Theorem 6.1 (1)]{U2}, and
so there exists a positive number $M$ (independent of $p$) such
that
$$|b(\nu_a, \nu_b)|<p \cdot |\FF| \cdot M p^{d+ \delta - \epsilon
-3} = M|\FF|p^{d+ \delta - \epsilon -2}.$$

\textbf{(Type II)} This is a node in $F_i \cap S_j$ with $j\neq
d+1$. Then again, we apply what we did in \cite[proof of Theorem
6.1]{U2}, to obtain the same estimate as above.

\textbf{(Type III)} This is a node in $F_i \cap S_{d+1}$. Since
$x_{d+1}=p-\sum_{k=1}^d x_k$, then we want $p-y_i'(p-\sum_{k=1}^d
x_k) \in \FF$ mod $p$, so $y_i'(\sum_{k=1}^d x_k) \in \FF$. But
this is again as in \cite[Theorem 6.1]{U2}, and we have the same
previous estimate. Notice that it works because $\delta-\epsilon
\geq 2$.

\textbf{(Type IV)} This is a node between $S_k$, $k \neq d+1$, and
a exceptional divisor over the fiber $F_j$. Notice that $\A_{p
\Delta}$ contains at least two fibers, so $\sum_{k=1}^d n_k x_k +
y_j < p$. Hence we are as in case (2) in the proof of
\cite[Theorem 6.1]{U2}.

\textbf{(Type V)} This is the new case, coming from nodes in the
resolution of singularities of $\A$. It does not involve
$x_{d+1}$. The idea is to analyze three equations $\EE_1$,
$\EE_2$, $\EE_3$ from the equation $\EE$, and estimate solutions
for each.

Without loss of generality, we rearrange indices so that $\nu_a=n
\sum_{i=1}^{\alpha} x_i + z$, for some $\alpha$, and $\nu_b=
\sum_{i=1}^{\alpha} x_i + \nu_a$, where $z =
\sum_{i=\alpha+1}^{\beta} n_i x_i + \varsigma y_j $, for some
$\beta$, with $\varsigma= 0$ or $1$, and $0<n_i<e$. Notice that $z
\neq 0$ for any solution of $\EE$. We define
$$\EE_1: \sum_{i=\beta+1}^d e x_i + \sum_{i \neq j} y_i =p_1,$$
equation with $m_1$ variables, $\EE_2: \sum_{i=1}^{\alpha} x_i =
p_2$ with $m_2$ variables, and $$\EE_3: \sum_{i=\alpha+1}^{\beta}
n_i x_i + \varsigma y_j  = p_3$$ with $m_3$ variables. So, $m_1 +
m_2 + m_3= d+\delta-\epsilon$. Notice that $p_i$ are numbers
varying in the region $0< p_i < p$, since we will look at
solutions of $\EE_i$ from solutions of $\EE$.

Say that $p-\nu_a'\nu_b \in \FF$, which means mod $p$, $(\sum_k
x_k)\nu_a' \in -\FF -1$. Of course the set $-\FF-1$ has same size
as $\FF$. We now use repeatedly the fact that the number of
nonnegative integer solutions of $a_1 z_1 + \ldots + a_m z_m=q$
for coprime $a_i$'s is $\frac{q^{m-1}}{(m-1)!a_1 a_2 \cdots a_m} +
O(q^{m-2})$ (see \cite{BGK01}). Let $p$ be large enough. Given
$0<p_3<p$, the number of solutions of $\EE_3$ is $< M_3
p^{m_3-1}$.

Now, the key observation is that mod $p$ we have
$$p_2(np_2+p_3)'\equiv \bar{p}_2(n\bar{p}_2+p_3)' \Rightarrow
p_2(n\bar{p}_2+p_3)\equiv \bar{p}_2(np_2+p_3) \Rightarrow p_2 p_3
\equiv \bar{p}_2 p_3 \Rightarrow p_2 \equiv \bar{p}_2$$ because
$p_3$ is not zero. In this way, we have to choose $p_2$ in a set
of size $|\FF|$. Now we fix $p_2$ and have at most $M_2 p^{m_2-1}$
solutions for $\EE_2$. After we have solutions for $\EE_3$ and
$\EE_2$, we have at most $M_1 p^{m_1 -1}$ solutions for $\EE_1$.
Putting it all together,
$$ b(\nu_a,\nu_b) < p \cdot M_3 p^{m_3-1} \cdot |\FF| \cdot M_2
p^{m_2-1} \cdot M_1 p^{m_1-1}= M_1 M_2 M_3 |\FF|
p^{d+\delta-\epsilon-2}.$$

But we know that $|\FF|< \sqrt{p} \big(\log(p)+2\log(2) \big)$
\cite[Apendix]{U2}, and that the total number of solutions of
$\EE$ is $\frac{p^{d+\delta-\epsilon-1}}{(d+\delta-\epsilon-1)!
e^d} + O(p^{d+\delta-\epsilon-2})$. Then, since the number of
nodes of $D_{\text{red}}$ is of course independent of $p$, we have
proved the existence of good solutions, and that a random one is
good with probability tending to $1$ as $p$ becomes arbitrarily
large.

Now, given good solutions with $p$ large, we proceed as at the end
of the proof of \cite[Theorem 6.1]{U2}, showing that \footnote{In
\cite{U2} there is a minor error for the estimate of $s(q,p)$.
This is due to the usual normalization by $12$ of a Dedekind sum.
The correct estimate is $12 |s(q,p)| \leq 3 \sqrt{p} + 5$, which
of course does not affect any asymptotic result.}
$$ LCF < \Big(\sum_{i<j} D_i.D_j \Big)(3\sqrt{p} + 2) \ \ \ \
\big| CCF \big| < \Big(\sum_{i<j} D_i.D_j \Big)(6\sqrt{p} + 7).$$
This proves the asymptotic result. Finally, these surfaces are of
general type because of the classification of algebraic surfaces
(see \cite{Beau96} for any characteristic), since we know that
$\bar{c}_1^2>0$ and $\bar{c}_2>0$ (Corollaries \ref{c61} and
\ref{c62}).
\end{proof}

\begin{obs}
With this theorem, one recovers the log Miyaoka-Yau inequalities
in Corollaries \ref{c61} and \ref{c62}, when $\K=\C$. We just
apply the (projective) Miyaoka-Yau inequality to the surfaces $X$
for large primes $p$. \label{o71}
\end{obs}

A good looking corollary, consequence of Theorem \ref{t41},
Section \ref{s6}, and Theorem \ref{t71}.

\begin{cor}
Let $d\geq 3$ be an integer, and let $A$ be an irreducible
projective curve in $\P^{d-2}(\K)$ not contained in the hyperplane
arrangement $\H_d$ (see above Example \ref{e21}). Then, there
exist nonsingular projective surfaces $X$ associated to $A$ such
that $X$ is of general type and $0<2c_2(X)<c_1^2(X)$.  \label{c71}
\end{cor}

\begin{proof}
Consider the arrangement $\A$ defined by $A$ as in Corollary
\ref{c41}. Then use Theorem \ref{t71} for $\A_{\Delta}$, and
Corollary \ref{c61}.
\end{proof}

\begin{obs}
In addition, one can prove that $\pi_1^{\text{\'et}}(X) \simeq
\pi_1^{\text{\'et}}( \overline{A})$, where $\overline{A}$ is the
normalization of the curve $A$, and $\pi_1^{\text{\'et}}$ denotes
the \'etale fundamental group \cite{U3}. \label{o72}
\end{obs}

\begin{cor}
Assume $\K=\C$. Let $\A$ be an arrangement of sections of $\pi:
\P_C(\L) \rightarrow C$. Then, there exist nonsingular projective
surfaces $X$ of general type such that $$2 <
\frac{c_1^2(X)}{c_2(X)} < 3, $$ (and so of positive index) having
$\frac{c_1^2(X)}{c_2(X)}$ arbitrarily close to
$\frac{\bar{c}_1^2(\A_{\Delta})}{\bar{c}_2(\A_{\Delta})}$. In
addition, there is an induced connected fibration $\pi': X
\rightarrow C$ which gives an isomorphism: $\pi_1(X) \simeq
\pi_1(C)$. In this way, Alb$(X) \simeq$ Jac$(C)$ and $\pi'$ is the
Albanese fibration of $X$. \label{c72}
\end{cor}

\begin{proof}
The first part is implied by Theorem \ref{t71} for $\A_{\Delta}$
and Corollary \ref{c61}. The map $f: X \rightarrow Y$ in the proof
of Theorem \ref{t71} is totally ramified along
$\bar{\A}_{\Delta}$, and $\pi \circ \sigma: Y \rightarrow C$ is a
connected fibration with at least one simply connected fiber and
one section in $\bar{\A}_{\Delta}$. Therefore, the construction
induces a connected fibration $\pi': X \rightarrow C$, and by
\cite[Proposition 8.3]{U2} we have $\pi_1(X) \simeq \pi_1(C)$. The
last part is a simple consequence of Albanese maps which applies
to any such fibration (see \cite{Beau96} for example).
\end{proof}

\begin{obs}
Corollary \ref{c72} is also valid for any $\A_{p\Delta}$ except
for $2<\frac{c_1^2(X)}{c_2(X)}$. If one thinks that the closest
$\frac{c_1^2(X)}{c_2(X)}$ is to $3$, the more interesting are the
surfaces $X$, then one may consider the construction starting with
some $\A_{p\Delta}$ (this is with $\epsilon>0$). For line
arrangements, this is indeed the case (see Remark \ref{o64}). By
Proposition \ref{p63}, we only need to consider primitive
arrangements in order to find upper bounds for Chern ratios.
\label{o73}
\end{obs}

\begin{example}
The conic in Example \ref{e51} produces an arrangement $\A \in
\AA(\P^1,\O(3),11)$. In the table of Example \ref{e61}, we
computed log Chern ratios for the extended and some partially
extend arrangements induced by $\A$. Then, by Corollary \ref{c72},
there are simply connected nonsingular projective surfaces of
general type with Chern ratios arbitrarily close to the ones in
that table. Notice that the highest is attained by a partially
extended arrangement, which avoids ``too many" double points.
\label{e71}
\end{example}

\begin{example}
Assume $\K$ has positive characteristic $p$. Take any $\A' \in
\AA(C,\L,d)$ for some curve $C$ and line bundle $\L$. Consider the
$\K$-linear Frobenius pull-back of $\A'$ composed $r$ times, as in
Remark \ref{o65}. Denote the resulting arrangement by $\A$. Then,
by Remark \ref{o65} and Theorem \ref{t71}, there are nonsingular
projective surfaces of general type $X$ with
$\frac{c_1^2(X)}{c_2(X)}$ arbitrarily close to
$\frac{\bar{c}_1^2(\A_{\Delta})}{\bar{c}_2(\A_{\Delta})}= 2 + p^r
\big(\frac{\bar{c}_1^2(\A'_{\Delta})}{\bar{c}_2(\A'_{\Delta})} -2
\big)$, and so arbitrarily large. We can prove that
$\pi_1^{\text{\'et}}(X) \simeq \pi_1^{\text{\'et}}(C)$ (see
\cite{U3}). Therefore, for any given positive characteristic and
nonsingular projective curve $C$, there are nonsingular projective
surfaces of general type $X$ with $\pi_1^{\text{\'et}}(X) \simeq
\pi_1^{\text{\'et}}(C)$ and violating any sort of Miyaoka-Yau
inequality. \label{e72}
\end{example}

\section{Appendix: log inequalities} \label{ap}

In this section, the ground field is $\C$. After fixing
$\AA(C,\L,d)$, it is clearly of our interest to find optimal upper
bounds for $\frac{\bar{c}_1^2}{\bar{c}_2}$ for extended and
partially extended arrangements (see Remark \ref{o64}).
Arrangements attaining upper bounds should be very special, and
they would produce interesting surfaces via Theorem \ref{t71}.

In this appendix, we show through Theorems \ref{apt} and
\ref{aptt} how this question about sharp upper bounds is connected
to old questions by Lang and others on effective height
inequalities \cite[pp.149-153]{Lang91}, via an inequality of Liu
\cite[Theorem 0.1]{Liu96}. Also, in a more general setting, we
show a way to obtain strictness for the log inequalities in
Corollaries \ref{c61} and \ref{c62}, and Theorems \ref{apt} and
\ref{aptt}. The next lemma follows from Kobayashi \cite{Ko85} and
Mok \cite{Mok09}.

\begin{lem}
Let $Y$ be a smooth projective surface, and let $D$ be a simple
normal crossing divisor in $Y$. Assume $K_Y + D$ is big and nef,
and $\bar{c}_1^2(Y,D)=3\bar{c}_2(Y,D)$.

Then, $D$ is a disjoint union of smooth elliptic curves.
\label{apl}
\end{lem}

\begin{proof}
By \cite[p.46]{Ko85}, $K_Y+D$ big and nef and
$\bar{c}_1^2(Y,D)=3\bar{c}_2(Y,D)$ imply that the universal
covering of $Y \setminus D$ is the complex two dimensional ball
$\B=\{(z,w) \in \C^2 : |z|^2+|w|^2 <1 \}$. Hence, there exist a
discrete group $\Gamma$ in $\aut (\B)$ such that $\B / \Gamma
\simeq Y \setminus D$. In particular, $\B / \Gamma$ has finite
volume. Notice that $\Gamma$ is torsion-free since it acts freely
on $\B$. Therefore, by \cite[Main Theorem ]{Mok09}, there exists a
smooth projective Mumford compactification $W$ of $\B / \Gamma$
such that $W \setminus (\B / \Gamma)$ is a disjoint union of
smooth elliptic curves $E_i$. In this way, we obtain a birational
map $W \dashrightarrow Y$. We now resolve this map and get a
birational morphism $\sigma: \widetilde{W} \rightarrow Y$. Then,
the inverse image $\widetilde{E}_i$ of each $E_i$ under $\sigma$
dominates $D_i$, after reordering indices. It is easy to see that
$\widetilde{E}_i=D_i$. But $\widetilde{E}_i$ is a smooth elliptic
curve with some finite trees of $\P^1$'s attached. Given one of
these trees, one has a smooth rational curve $F$ intersecting
$\widetilde{E}_i-F$ at one point. But then $0 \leq (K_Y + D).F=-2
+ 1 = -1$. So, there are no trees, and $E_i=D_i$ for all $i$.
\end{proof}

Let $f: Y \rightarrow C$ be a fibration of a smooth projective
surface over a smooth projective curve $C$, denote by $g$ the
genus of the generic (connected) fiber of $f$ and by $q$ the genus
of $C$. Let $\omega_{Y|C}:=K_Y - f^*(\omega_C)$ be the relative
dualizing sheaf.

Let $S_1,\ldots,S_n$ be $n$ mutually disjoint sections of $f$.
Assume $f$ is a semi-stable fibration of $n$-pointed curves of
genus $g$, marked by these sections. Let
$$D=S_1+ \ldots +S_n + f^*(c_1+\ldots+c_{\delta})$$ where $c_1,\ldots,c_{\delta}$ are the images of the singular fibers of $f$.

\begin{teo}
Let $f$ be not isotrivial, i.e., the moduli of its fibers varies
as $n$-pointed semi-stable curves. Assume $D\neq \emptyset$, and
$n\geq 1$ when $g=1$. Then $$0 < \bar{c}_1^2(Y,D) < 3
\bar{c}_2(Y,D).$$ \label{apt}
\end{teo}

\begin{proof}
The generic fiber has $\bar{\kappa}\big(f^{-1}(c), (S_1+ \ldots
+S_n )|_{f^{-1}(c)} \big)=1$ ($\bar{\kappa}$ denotes the log
Kodaira dimension): $\P^1$ minus at least four points or elliptic
curve minus at least one point or the rest. Now, since $f$ is not
isotrivial, it has at least $3$ singular fibers when $C=\P^1$ (see
\cite{Beau81}), or at least one when $C$ is an elliptic curve (see
\cite[p.127]{BHPV04}). So, in any case, the base is of log general
type (on the base we take the log curve
$(C,c_1+\ldots+c_{\delta})$). By a theorem of Kawamata
\cite[Theorem 11.15]{IitakaAG82}, we have for a general $c \in C$
$$ \bar{\kappa}(Y,D) \geq \bar{\kappa}\big(f^{-1}(c), (S_1+ \ldots
+S_n )|_{f^{-1}(c)} \big) +
\bar{\kappa}(C,c_1+\ldots+c_{\delta}),$$ and so $(Y,D)$ is of log
general type. Notice that $D$ is a semi-stable curve, just because
the fibration is semi-stable, and when $C=\P^1$ we have at least
$3$ singular fibres. Therefore, Sakai's theorem \cite{Sakai1}
applies, and so $\bar{c}_1^2(Y,D) \leq 3 \bar{c}_2(Y,D)$. Below we
show that $K_Y+D$ is nef to obtain the strict inequality.

(case $g=0$): We can explicitly show that $K_Y+D$ is nef
(Corollary \ref{c61}).

(case $g\geq 2$): Let $g \colon Y \rightarrow Y'$ be the relative
minimal model of $f$. Let $E_i$ be the exceptional divisors. We
can write $$K_Y +D \equiv g^*(\omega_{Y'|C}) + \sum_k n_k E_k +
f^*(\omega_C) + D$$ where for some positive integers $n_i$'s.
Since $g\geq 2$, the dualizing sheaf $\omega_{Y'|C}$ is nef (due
to Arakelov). If $q>0$, then $f^*(\omega_C)$ is nef as well, and
so we check that $K_Y +D$ is nef by intersecting it with
components of $D$ (notice that $D$ includes $E_k$'s). If $q=0$, we
have at least $3$ singular fibers, and so we delete
$f^*(\omega_C)$ using $D$, and again to check nef we intersect
$K_Y +D$ with the components of $D$.

(case $g=1$): As in the previous case, we go to $Y'$. Notice that
$Y'$ has no multiple fiber. By the canonical bundle formula we
have \cite[p.214]{BHPV04} $$K_Y +D \equiv (\chi(Y)+2(q-1))F +
\sum_k n_k E_k + D$$ where $F$ is a general fiber of $f$. So, if
$q>0$, we are done by the previous argument. If $q=0$, we are done
by the same argument, since there are at least $3$ singular fibers
by \cite{Beau81}.

Therefore, $K_Y + D$ is nef, and strict inequality follows from
Lemma \ref{apl}.
\end{proof}

Theorem \ref{apt} is also valid when $D= \emptyset$ (i.e., a
Kodaira fibration). It follows from \cite[Theorem 0.1]{Liu96}. Our
argument does not prove it. Actually, up to the case
$D=\emptyset$, Theorem \ref{apt} is just a small extension of
\cite[Theorem 0.1]{Liu96} (since we also consider the cases
$g=0,1$), as we now see. We have $\bar{c}_2(Y,D) = e(X) - e(D)$
($e(A)$ is the Euler topological characteristic of $A$) and the
usual formula \cite[Lemma VI.4]{Beau96} $$e(Y)= 4(g-1)(q-1) +
\sum_{i=1}^{\delta} \big( e(f^{-1}(c_i)) - e(F) \big) $$ where $F$
is a generic fiber of $f$. One sees that $e(D)=
\sum_{i=1}^{\delta} e(f^{-1}(c_i)) +n(2-2q) - \delta n$. So,
$\bar{c}_2(Y,D) = (2g-2+n)(2q-2+\delta).$ Obviously $\omega_{Y|C}
+ D = K_Y + D - f^* \big(\omega_C + \sum_{i=1}^{\delta} c_i
\big)$. We then square it and see the inequality in \cite[Theorem
0.1]{Liu96}.

\begin{obs} In \cite[Theorem 0.1]{Liu96}, where $g \geq 2$ is assumed,
we have that $f$ is isotrivial if and only if $\bar{c}_1^2(Y,D)=3
\bar{c}_2(Y,D)$. But we now show that this corresponds to
uninteresting situations. First notice that $K_Y +D$ is nef by the
same argument used in (case $g\geq 2$) of Theorem \ref{apt}. If
$\bar{c}_1^2(Y,D)>0$, then $K_Y +D$ becomes big and nef, and we
apply Lemma \ref{apl} to obtain a contradiction, unless $q=1$ and
$\delta=0$. But then $\bar{c}_2(Y,D)=0$, which is a contradiction
to our assumption $\bar{c}_1^2(Y,D)>0$. Therefore, we are in the
trivial case $\bar{c}_1^2(Y,D)= \bar{c}_2(Y,D)=0$. \label{apo1}
\end{obs}

For completeness' sake, we explicitly show the connection with
height inequalities of algebraic points on curves over function
fields. This is another proof of Tan's height inequality
\cite[Theorem A]{Tan95}. Let $f: Y \rightarrow C$ be a connected
fibration as before, denoting by $g$ the genus of the generic
fiber of $f$ and by $q$ the genus of $C$. Assume that $f$ is
semi-stable. Let $\text{K}(C)$ be the function field of $C$. For
an algebraic point $P \in Y(\overline{\text{K}(C)})$, let $C_P$ be
the corresponding horizontal curve (i.e. multisection) in $Y$. As
usual, let
$$h_K(P)= \frac{\omega_{Y|C}.C_P}{F.C_P} \ \ \ \ \ \ \
d(P)= \frac{2g(\overline{C_P})-2}{F.C_P} $$ be the geometric
height and the geometric logarithmic discriminant respectively.
The curve $\overline{C_P}$ is the normalization of $C_P$, and $F$
is a general fiber of $f$.

\begin{teo}
Assume $g\geq 2$, and that $f$ is not isotrivial. Let $\delta$ be
the number of singular fibers of $f$. Then, for any algebraic
point $P$, we have
\begin{center} $h_K(P)< (2g-1)(d(P)+\delta) - \omega_{Y|C}^2 .$ \end{center} \label{aptt}
\end{teo}

\begin{proof}
Let $C_P$ be the horizontal curve in $Y$ defined by $P$, and let
$g: \overline{C_P} \rightarrow C$ be the composition of the
normalization of $C_P$ with $f$, so $d:=\deg(g)=F.C_P$. Then, we
have
\begin{center} $ \xymatrix{ Y_P \ar[rd]^{f_P} \ar[r] & \overline{Y} \ar[d]^{\bar{f}} \ar[r]^{G} & Y \ar[d]^{f}   \\ & \overline{C_P}
\ar[r]^{g} & C }$
\end{center} where $\bar{f}$ is the unique semi-stable fibration induced by $g$. Notice that $G^*(C_P)$ contains a section $S$ of $\bar{f}$, by
construction. The map $f_P$ is the induced semi-stable fibration
with a marked point (marked by $S$). Let $\delta_P$ be the number
of singular fibers of $f_P$. Notice that $\delta_P$ is at most $d
\delta$. Consider $D=S' + f_P^*(c_1+\ldots+c_{\delta_P})$ where
$c_1,\ldots,c_{\delta_P}$ are the images of the singular fibers of
$f_P$ in $\overline{C_P}$, and $S'$ is the strict transform of
$S$. We now apply Theorem \ref{apt} to have $(K_{Y_P} + D)^2 < 3
(2g-2+1)(2g(\overline{C_P})-2 +\delta_P)$. But, one checks that
$(K_{Y_P}+D)^2=(K_{\overline{Y}} + S +
\bar{f}^*(c_1+\ldots+c_{\delta_P}))^2=(\omega_{\overline{Y}|{\overline{C_P}}}
+ S + (\delta_P + 2g(\overline{C_P})-2)F)^2$. Also, since
$\bar{f}$ is semi-stable, we know that
$G^*(\omega_{Y|C})=\omega_{\overline{Y}|{\overline{C_P}}}$, and by
the projection formula
$S.\omega_{\overline{Y}|{\overline{C_P}}}=C_P.\omega_{Y|C}$. So,
the log inequality above becomes \begin{center} $d \omega_{Y|C}^2
+ \omega_{Y|C}.C_P + 2(2g-1)(2g(\overline{C_P})-2+ \delta_P)  < 3
(2g-1)(2g(\overline{C_P})-2 + \delta_P)$,
\end{center} and so we rearrange to obtain the claimed height inequality (also use $\delta_P \leq d \delta$).
\end{proof}


\vspace{0.1 cm} {\small Department of Mathematics and Statistics,
University of Massachusetts at Amherst, USA.}

\begin{thebibliography}{99}

\bibitem{BHPV04}
    W. P. Barth, K. Hulek, C. A. M. Peters, and A. Van de Ven.
    \emph{Compact complex surfaces},
    Ergebnisse der Mathematik und ihrer Grenzgebiete. 3. Folge., second edition, vol. 4, Springer-Verlag, Berlin, 2004.

\bibitem{Beau81}
    A. Beauville.
    \emph{Le nombre minimum de fibres singuli\`{e}res d'une courbe sur $\P^1$},
    Ast\'{e}risque 86, 97--108(1981).

\bibitem{Beau96}
    A. Beauville.
    \emph{Complex algebraic surfaces},
    London Mathematical Society Student Texts, vol. 34, Cambridge University Press, Cambridge, 1996.

\bibitem{BGK01}
    M. Beck, I. M. Gessel, and T. Komatsu.
    \emph{The polynomial part of a restricted partition function related to the {F}robenius problem},
    Electron. J. Combin. 8 (2001), no. 1, Note 7, 5 pp. (electronic).

\bibitem{BruErd1948}
    N. G. de Bruijn, and P. Erd$\ddot{o}$s.
    \emph{On a combinatorial problem},
    Nederl. Akad. Wetensch., Proc. 51, (1948) 1277--1279 $=$ Indagationes Math. 10, 421--423 (1948).

\bibitem{CT1}
    A.-M. Castravet and J. Tevelev.
    \emph{Exceptional loci on $\overline{M}_{0,n}$ and hypergraph curves},
    arXiv:0809.1699v1.

\bibitem{DM}
    Deligne P. and Mumford D.
    \emph{The irreducibility of the space of curves of given genus},
    Inst. Hautes \'Etudes Sci. Publ. Math. 36 (1969), 75--109.

\bibitem{Ha1}
    R. Hartshorne.
    \emph{Algebraic geometry},
    Graduate Text in Mathematics v.52, Springer, 1977.

\bibitem{Hi1}
    F. Hirzebruch.
    \emph{Arrangements of lines and algebraic surfaces},
    Arithmetic and geometry, Vol. II, Progr. Math. 36, Birkh\"auser, Boston, Mass., 1983, 113--140.

\bibitem{IitakaAG82}
    S. Iitaka.
    \emph{Algebraic Geometry},
    Graduate Texts in Mathematics, vol. 76, Springer-Verlag, New York, 1982.

\bibitem{Ka1}
    M. M. Kapranov.
    \emph{Veronese curves and Grothendieck-Knudsen moduli space $\overline{M_{0,n}}$},
    J. Algebraic Geom. 2 (1993), 239--262.

\bibitem{Ka2}
    M. M. Kapranov.
    \emph{Chow quotients of Grassmannians I},
    Adv. Soviet Math. 16, part 2, A.M.S., (1993) 29--110.

\bibitem{Knudsen}
     F. F. Knudsen.
     \emph{The projectivety of moduli spaces of stable curves II: the stacks $M_{g,n}$},
     Math. Scand. 52 (1983), 161--199.

\bibitem{Ko85}
     R. Kobayashi.
     \emph{Einstein-{K}aehler metrics on open algebraic surfaces of general type},
     Tohoku Math. J. (2), Vol. 37, No.1 (1985), 43--77.

\bibitem{Lang91}
     S. Lang.
     \emph{Number Theory III: diophantine geometry},
     Encyclopaedia of Mathematical Sciences: V. 60 (1991).

\bibitem{Liu96}
     K. Liu.
     \emph{Geometric height inequalities},
     Math. Res. Lett. 3, 693--702 (1996).

\bibitem{Mok09}
     N. Mok.
     \emph{Projective-algebraicity of minimal compactifications of complex-hyperbolic space forms of finite volume},
     pre-print in http://hkumath.hku.hk/$\sim$imr/IMRPreprintSeries/2009/IMR2099-4.pdf (2009).

\bibitem{Sakai1}
    F. Sakai.
    \emph{Semi-stable curves on algebraic surfaces and logarithmic pluricanonical maps},
    Math. Ann. 254(1980), no. 2, 89--120.

\bibitem{Tan95}
    S.-L. Tan.
    \emph{Height inequality of algebraic points on curves over functional fields},
    J. Reine Angew. Math. 461(1995), 123--135.

\bibitem{U}
    G. Urz\'ua.
    \emph{Arrangements of curves and algebraic surfaces},
    Ph.D. Thesis, University of Michigan (2008).

\bibitem{U1}
    G. Urz\'ua.
    \emph{On line arrangements with applications to $3$-nets},
    Adv. Geom. 10(2010), 287--310.

\bibitem{U2}
    G. Urz\'ua.
    \emph{Arrangements of curves and algebraic surfaces},
    J. of Algebraic Geom. 19(2010), 335--365.

\bibitem{U3}
    G. Urz\'ua.
    \emph{Simply connected surfaces in positive characteristic},
    in preparation (2009).

\end{thebibliography}
\end{document}